\pdfoutput=1
\documentclass{article}
\usepackage[utf8]{inputenc}
\usepackage[T1]{fontenc}
\usepackage[english]{babel}
\usepackage[top=2.5cm, bottom=2.54cm, left=2.54cm, right=2.54cm]{geometry}
\usepackage{amsmath}
\usepackage{amssymb}
\usepackage{array}
\usepackage{listings}
\usepackage{framed}
\usepackage[framed, amsmath, thmmarks, amsthm]{ntheorem}
\usepackage{graphicx}
\graphicspath{{./}{figures/}}
\usepackage{hyperref}
\usepackage{caption}
\usepackage{fancyhdr}
\usepackage{multicol}
\usepackage{appendix}
\usepackage{verbatim}
\usepackage{tikz}
\usepackage{subfig}
\usepackage{float}
\usepackage{algorithm}
\usepackage{algpseudocode}
\def\restriction#1#2{\mathchoice
              {\setbox1\hbox{${\displaystyle #1}_{\scriptstyle #2}$}
              \restrictionaux{#1}{#2}}
              {\setbox1\hbox{${\textstyle #1}_{\scriptstyle #2}$}
              \restrictionaux{#1}{#2}}
              {\setbox1\hbox{${\scriptstyle #1}_{\scriptscriptstyle #2}$}
              \restrictionaux{#1}{#2}}
              {\setbox1\hbox{${\scriptscriptstyle #1}_{\scriptscriptstyle #2}$}
              \restrictionaux{#1}{#2}}}
\def\restrictionaux#1#2{{#1\,\smash{\vrule height .8\ht1 depth .85\dp1}}_{\,#2}} 

\newcommand{\fonction}[5]{\begin{array}[t]{lrcl}#1 :&#2 &\longrightarrow &#3\\&#4& \longmapsto &#5 \end{array}}

\newcommand{\bC}{\mathbb{C}}

\newcommand{\bE}{\mathbb{E}}

\newcommand{\bP}{\mathbb{P}}
\newcommand{\bR}{\mathbb{R}}

\newcommand{\D}{\mathcal{D}}
\newcommand{\Od}{\mathcal{O}}
\newcommand{\dom}{\D \setminus \Od}

\newcommand{\C}{\mathcal{C}}

\newcommand{\Tau}{\mathcal{T}}

\newcommand{\HH}{\mathrm{H}}
\newcommand{\WW}{\mathrm{W}}
\newcommand{\LL}{\mathrm{L}}
\newcommand{\DD}{\mathrm{D}}
\newcommand{\RR}{\mathrm{R}}
\newcommand{\nn}{\mathrm{n}}
\theoremstyle{break}
\newtheorem{theorem}{Theorem}[section]
\newtheorem{lemme}[theorem]{Lemma}

\newtheorem{remark}[theorem]{Remark}
\newtheorem{definition}[theorem]{Definition}
\newtheorem{proposition}[theorem]{Proposition}

\numberwithin{equation}{section}

\usepackage{authblk}
\usepackage{subfloat}
\usetikzlibrary{shapes.geometric, arrows.meta, positioning}

\tikzstyle{step} = [rectangle, draw, rounded corners, minimum width=5cm, minimum height=1cm, text centered, fill=blue!10]
\tikzstyle{loop} = [rectangle, draw, rounded corners, minimum width=6cm, minimum height=1cm, text centered, fill=orange!20]
\tikzstyle{arrow} = [thick,->,>=stealth]

\title{Simultaneous reconstruction of an obstacle and its wavenumber for the Helmholtz equation: a robust optimization approach}

\author[1]{Thomas Bonnafont\thanks{thomas.bonnafont.enac@gmail.com}}
\author[2]{Fabien Caubet\thanks{fabien.caubet@univ-pau.fr}}
\author[2]{Nicolas Grima\thanks{nicolas.grima@univ-pau.fr}}
\author[2]{Jason Tridon\thanks{jason.tridon@univ-pau.fr}}

\affil[1]{\normalsize Lab-STICC, UMR CNRS 6285, ENSTA, Institut Polytechnique de Paris,  Brest, France}
\affil[2]{\normalsize Universit\'e de Pau et des Pays de l'Adour, E2S UPPA, CNRS, LMAP, UMR 5142, 64000 Pau, France.}

\date{\today}

\begin{document}

\maketitle

\begin{abstract}
In this paper, we propose a method to recover simultaneously the shape and wavenumber of an object from noisy radar or sonar measurements. The physical problem is modeled in a bounded domain with the Helmholtz equation and appropriate transmission conditions. The approach relies on a Kohn–Vogelius formulation of the inverse problem and uses Nesterov’s accelerated gradient descent combined with shape optimization techniques. To address the robustness of the method with respect to measurement noise, we are looking to minimize a convex combination of the expected value and variance, which we achieve by using a Karhunen–Loève expansion. The effectiveness of the approach is demonstrated through two-dimensional numerical experiments.
\end{abstract}

\bigskip


\medskip


\section{Introduction}

\subsection{Motivations}
The inverse scattering problem~\cite{cakoni2014qualitative} consists in retrieving the physical and parametric properties of an object from radar or sonar measurements. This problem is of major interest in numerous applications, including defense, non-destructive testing, health monitoring, and geoscience. For a general mathematical background on this problem, we refer to the works of Colton and Kress~\cite{colton2013,cakoni2014qualitative,coltonkress2019,cakoni2022inverse}. \\

This work is situated within this context, with the objective of detecting  and reconstructing both the shape and the wavenumber of an unknown object from radar measurements while accounting for uncertainties. We expect this strategy to overcome the problem of ill-posedness of retrieving the wavenumber alone~\cite{colton2013,cakoni2014qualitative}. Here, as usual, we consider the case of monochromatic incident waves and thus deal with the Helmholtz equation, with Robin boundary condition. This outer boundary condition models the Sommerfeld one on a bounded domain. Finally, since we are interested in reconstructing the shape of the target and its wavenumber, the governing equations involve transmission conditions on the unknown interface between the object and the exterior environment. Our transmission condition models a situation in which the dielectric permittivity $\varepsilon$ of the object differs from that of the surrounding medium. \\

The adopted strategy is the minimization of an objective functional of Kohn-Vogelius type using Nesterov's accelerated gradient descent. From a numerical point of view, we parameterize the shape by a truncated Fourier series. This allows us to work naturally in a finite-dimensional space. Our methodology thus relies on tools from shape optimization, presented in the books of Allaire~\cite{allaire2007}, and Henrot and Pierre~\cite{henrotpierre2018}.  \\

Previous studies have considered the shape reconstruction for a similar problem. We can highlight the work~\cite{Litman1998} where an exterior Sommerfeld radiation condition is employed. As previously mentioned, we approximate this condition by the Robin condition. This specific choice allows us to efficiently solve the underlying equations using a finite element method. A further key difference is that they, along with many other authors (e.g., \cite{Dorn2000, Audibert2022}), employ level-set methods, that allow one to consider topological change. However, this does not address the ill-posedness of the inverse scattering problem while our parametrization addresses this issue naturally. \\

Another strand of the literature is the linear sampling method (see, e.g.,~\cite{bucci1997electromagnetic,cakoni2011linear,crocco2012linear}). The method is mostly based on a linearization of the integral equation in order to invert the system and retrieve the indicator of its shape, even with unknown dielectric parameters. Nonetheless, this method is still limited and the results highly depend on the contrast between the target and the surrounding domain. Finally, we highlight that this approach does not address the reconstruction of the wavenumber.\\

Here, we do not consider thin-layers, but in another context, the transmission condition can be replaced by modeling it through a Ventcel-type condition~\cite{kress2018,Boutarene2024}. Various studies have investigated the simultaneous reconstruction of the thickness and the impedance condition of the unknown object as~\cite{Bourgeois2012}. Those analysis does not incorporate data uncertainties into the reconstruction model. Also, we can remark that in~\cite{Bourgeois2012}, the authors prove the following identifiability result: the solution to the inverse problem is unique assuming that any directions of observations and scattering are known. Similarly, we refer to the works~\cite{Nachman}, \cite{Novikov}, \cite{hahner}, \cite{coltonkress2019} and \cite{sylvesteruhlmann} for useful papers about uniqueness in inverse scattering problem and its link with electrical impedance tomography. \\

The use of the Kohn-Vogelius functional has previously been considered in other physical problems, see~\cite{Caubet2013,Dambrine2022,Chaabane2023}. Indeed, this approach has shown to be better than the conventional least-square one, even for a simple inverse scattering case where homogeneous Dirichlet condition is considered for the target boundary~\cite{bonnafont2023inverse}. In particular, in~\cite{Dambrine2022} they introduced a Karhunen–Loève expansion to handle uncertainties in the measurements, which we similarly employ in the present study. The authors also seek to reconstruct an object using shape optimization. However, their work is set within the framework of obstacle reconstruction with Poisson equation.\\

This paper is organized as follows. Section~\ref{sec:mainres} presents the equations under study and the main results we obtained. We introduce a Kohn-Vogelius, which is specifically constructed so that its minimizers coincide with the solutions of the inverse problem. We give an identifiability result in dimension~$3$ and prove the existence and give the expressions of the derivatives for both the deterministic case and the case with uncertainties. In the latter case, we develop the Karhunen–Loève expansion that permits us to replace the cost functional by moments of order $1$ and $2$ of the deterministic Kohn-Vogelius functional. Section~\ref{sec:numexp} is devoted to illustrate our method by numerical experiments. In particular, we highlight the advantages of our method. In these experiments, the objects are represented by truncated Fourier expansions and we minimize the cost functional, thanks to a Nesterov descent scheme. Section~\ref{sec:proofs} shows that the equations considered are well-posed and addresses the proofs of the principal results. 

\subsection{Notations and setting of the inverse problem}

\paragraph{Notations.}
As previously mentioned, in this paper, we focus on an inverse scattering problem. In particular, even if the work is general, we use the notation of electromagnetic. The situation can be stated quite simply: we want to detect and reconstruct an unknown obstacle $\Od$ (which could be, for instance, a tumor or an airplane), assumed to be included in a domain $\D$, and its wavenumber, from electromagnetic measurements. \\

\begin{figure}[h!tp]
\centering
\begin{tikzpicture} 

\filldraw[fill=gray!20, draw=gray!20] plot [smooth cycle, tension=1] coordinates {(0,0) (0.8,0.05) (2,1) (2.5,2) (3,3) (2.25,3.5) (1.5,4) (0.75,4) (0,4) (-0.75,4) (-1.5,4) (-1.75,2.5) (-2,1) (-1,0.) (0,0)};

\filldraw[fill=blue!20, draw=blue!20] plot [smooth cycle, tension=1] coordinates {(0.5,1.5) (1,1.75) (1.5,2) (1.35,2.5) (1.2,3) (0.7,3) (0.2,3) (-0.15,2.75) (-0.5,2.5) (0,2)};

\begin{scope}
\clip (0,-0.3) -- (1.5,4.5) -- (3.1,4.5) -- (3.1,-0.3); 
\draw[orange, thick] plot [smooth cycle, tension=1] coordinates {(0,0) (0.8,0.05) (2,1) (2.5,2) (3,3) (2.25,3.5) (1.5,4) (0.75,4) (0,4) (-0.75,4) (-1.5,4) (-1.75,2.5) (-2,1) (-1,0.) (0,0)};
\end{scope}
\begin{scope}
\clip (0,-0.3) -- (1.5,4.5) -- (-2.1,4.5) -- (-2.1,-0.3); 
\draw[blue, thick] plot [smooth cycle, tension=1] coordinates {(0,0) (0.8,0.05) (2,1) (2.5,2) (3,3) (2.25,3.5) (1.5,4) (0.75,4) (0,4) (-0.75,4) (-1.5,4) (-1.75,2.5) (-2,1) (-1,0.) (0,0)};
\end{scope}

\node[orange] at (2.35, 2.5) {$\Gamma_{\mathrm{m}}$};
\node[orange] at (2.8,1.5) {$u_{\mathrm{obs}}$};
\node[blue] at (-2., 3.5) {$\Gamma_{\mathrm{c}}$};

\node at (-1.2, 2.5) {$\mathcal{D}$};
\node at (0.5, 2.5) {$\mathcal{O}$};

\draw[->, thick] (3.4, 4.8) -- (2.4, 3.8);
\draw[->, thick] (3.6, 4.8) -- (2.6, 3.8);
\draw[->, thick] (3.2, 4.8) -- (2.2, 3.8) node[midway, left] {$u_{\mathrm{inc}}$};

\end{tikzpicture}
\caption{An example of scattering inverse problem.} 
\label{fig:enter-label}
\end{figure}
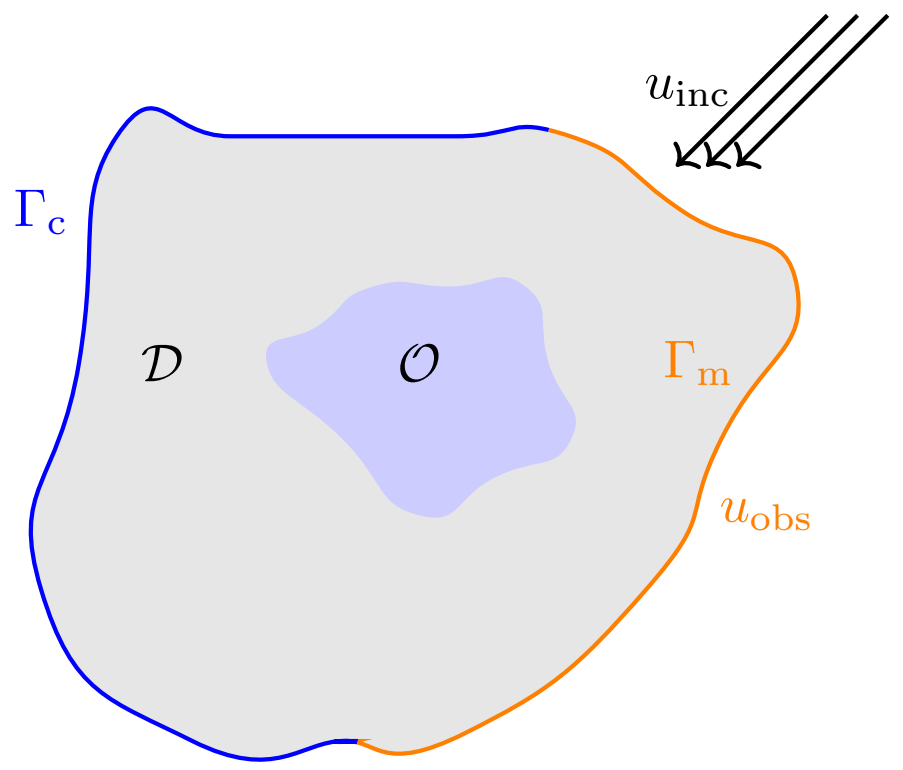

We study the harmonic case, and thus assume an $\exp(-i2\pi f t)$-time variation where $f$ is the frequency, inducing the $\exp(ikr)$ space variations. In all the following, $\D$ denotes a bounded domain, i.e. a connected open subset of $\bR^d$ $(d=2,3)$ with Lipschitz boundary. To detect the object $\Od$, we send planar waves on this object and we measure a signal $u_{\mathrm{obs}}$ via some sensors (such as antennas) on part of the boudary of $\D$ here denoted by~$\Gamma_{\mathrm{m}}$. The boundary $\Gamma_{\mathrm{c}}$ denotes the part of the boundary that is not covered by the sensors. Hence, we have $\partial \D = \overline{\Gamma_{\mathrm{m}} \cup \Gamma_\mathrm{c}}$. Both parts are not necessarily connected. The Robin condition is imposed on $\partial \D$. The wave sent corresponds to the incident wave $u_{\mathrm{inc}}$ and the wave that comes back to us is the sum of the scattered wave $u_{\mathrm{s}}$ and $u_{\mathrm{inc}}$, that is $u = u_{\mathrm{inc}} + u_\mathrm{s}$. Hence, we know $ u_{\mathrm{inc}}$ everywhere and we also know the sum $u = u_{\mathrm{obs}} = u_{\mathrm{inc}} + u_{\mathrm{s}}$ on the sensors $\Gamma_{\rm m}$. For a schematic representation of the inverse scattering problem, one can refer to Figure~\ref{fig:enter-label}. \\

In the following, we will consider the following set of admissible shapes $\Od$:
$$
\mathcal{U}_{\mathrm{adm}} = \left\lbrace \Od \subset \D ; \; \Od \mbox{ is closed with Lipschitz boundary and such that } \dom \mbox{ is connected} \right\rbrace .
$$ 

We also assume that the wavenumber is constant in both $\Od$ and $\D\subset \Od$. From now on, we write $\Od_1$ the unknown and $\Od_2 = \D \setminus \Od_1$ such that the wavenumber on $\Od_j$ is written $k_j$, $j=1,2$. Notice that $k_{1}$ is then unknown and $k_{2}$ is known.
In the following, we assume that
$$
\mathfrak{Im}(k_{2}) > 0
$$
and we will consider the following set of admissible wavenumber $\mathfrak{K}_{\mathrm{adm}}$:
$$
\mathfrak{K}_{\mathrm{adm}} = \left\lbrace k_{1} \in \mathbb{C} ; \; \mathfrak{Re}(k_{1})\mathfrak{Re}(k_{2}) > 0 \mbox{ and } \mathfrak{Im}(k_{1}) > 0 \right\rbrace  .
$$ The conditions on the imaginary parts of the wavenumbers are in line with the convention we choose for the time-variation. It shall also be noted that this also ensure the condition on the real part of $k_1$ and $k_2$ since they will be of same sign. In the other case, both imaginary parts are negative provided that we also change the sign of the term $-iku$ in the Robin condition. \\

\paragraph{Setting of the inverse problem.}
The inverse problem can then be stated as follows: 
we aim to find an admissible shape $\Od_1 \in \mathcal{U}_{\mathrm{adm}}$, an admissible wavenumber $k_1 \in \mathfrak{K}_{\mathrm{adm}}$ and $u \in V$ such that
$$
\left\{
\begin{array}{lll}
    \Delta u_j + k_j^2 u_j = 0 & \text{ in } \Od_j \\
    \partial_{\rm n} u_{2} - ik_2u_{2} = h & \text{ on } \partial \D \\
    u_{2} = u_{\mathrm{obs}} & \text{ on } \Gamma_{\mathrm{m}} \\
    \left[u\right] = [\partial_{\rm n} u] = 0 & \text{ on }  \partial \Od_1 ,
\end{array}
\right.
$$
where $h \in \HH^{-1/2}(\partial\D)$ and $u_{\mathrm{obs}} \in \HH^{1/2}(\Gamma_{\mathrm m})$ are given and 
$
V = \{u \in \mathrm{H}^1(\D) \text{ such that } [u] = 0 \text{ on } \partial \Od_1\}.
$
Here and in all the following, for each quantity $o_j$ defined on a subdomain $\Od_j$, we denote by $[o] = o_1 - o_2$ the jump of $o$ on $\partial \Od_1$. In order to simplify the notations, in the following, we also denote by $\chi_j$ the characteristic function of $\Od_i$ and omitting the subscript implies that we consider the sum $o = \chi_1 o_1 + \chi_2 o_2$. As an important example, $k$ denotes the sum $\chi_1 k_1 + \chi_2 k_2$ and we can write 
$$
\Delta u + k^2 u = 0 \; \mbox{ in } \D.
$$

\begin{remark}
We underline that, using the Robin boundary condition on $\partial\D$, both Dirichlet $u_{\mathrm{obs}}$ and Neumann $\partial_{\rm n}u_{\mathrm{obs}}$ conditions are known on $\Gamma_{\mathrm{m}}$.
\end{remark}

\section{Main results}
\label{sec:mainres}
In this section, we define the Kohn-Vogelius functional used in this work and state the main theoretical results of this work. First, we calculate its derivative with respect to its two unknowns, that are~$\Od_1$ and~$k_1$. Then, we give a method that takes into account noise on the input data, i.e., the robust method. The general forward problem we aim to solve can be stated as finding functions $u_j \in \HH^1(\Od_j)$ such that
\begin{equation}\tag{$\mathcal{T}(\Od_1)$}\label{eq:transmission}
\left\{ \begin{array}{rcll}
    \Delta u + k^2 u & = & f & \text{ in } \mathcal{D} \\
    \partial_{\nn} u - ik_2u & = & h & \text{ on } \Gamma_{\mathrm{c}} \\
    u & = & u_{\mathrm{obs}} & \text{ on } \Gamma_{\mathrm{m}} \\
    \left[u\right] & = & 0 & \text{ on } \partial \Od_1 \\
    \left[\partial_{\nn} u\right] & = & g & \text{ on }  \partial \Od_1 ,
    \end{array}
\right.\end{equation} 
where $g \in \HH^{-1/2}(\partial \Od_1), h \in \HH^{-1/2}(\Gamma_{\rm c}), f \in \mathrm{L}^{2}(\D)$ and $u_{\mathrm{obs}} \in \HH^{1/2}(\Gamma_{\mathrm m})$. The potential errors and uncertainties will occur on the functions $h$ and $u_{\rm obs}$. 

\begin{remark}
The physical forward case corresponds to $f = 0$, $g = 0$ with empty $\Gamma_{\rm m}$. This means that we do not consider any sensors or source terms and a transmission condition $\left[\partial_{\nn} u\right] = 0$, with an incident wave contained in~$h$. However, we give here a general version of the problem in order to show the well-posedness of all the equations written below (in particular the adjoint problems~\eqref{eq:adjoints}).
\end{remark}

The homogeneous forward problem (with $u_{\mathrm{obs}} = 0$) \eqref{eq:transmission} writes under variational form as finding $u \in V$ such that
\[
\int_{\D} \left( \nabla u \cdot \nabla \overline{v} - k^2 u \overline{v} \right) - ik_2 \int_{\partial \D} u \overline{v} = \int_{\partial \D} h \overline{v} + \int_{\partial \Od_1} g \overline{v} - \int_{\D} f \overline{v} \quad \text{ for all } v \in V,
\] where the notation $\overline{v}$ denotes the complex conjugate of $v$. Note that since our functions are complex-valued, the scalar product we consider is hermitian. 

\begin{proposition}\label{prop:wpforward}
Under the assumptions above, the transmission problem $(\mathcal{T}(\Od_1))$ is well-posed.
\end{proposition}

Concerning the above inverse problem, we now reformulate it as an optimization problem. We introduce the Kohn-Vogelius (also used in other example as in~\cite{Caubet2013}) functional as follows

\begin{equation}\label{def:K}
    \mathcal{K}(\Od_1,k_1) =  \frac{1}{2}\int_{\D} |\nabla(u_{\DD} - u_{\RR})|^2,
\end{equation}
where $u_{\mathrm{D}} = u_{\DD}(\mathcal{O}_1,k_1) \in V$ and $u_{\mathrm{R}} = u_{\RR}(\mathcal{O}_1,k_1) \in V$ are the unique solutions to the problems

\begin{equation} \label{eq:urud}\begin{array}{cc}
\left\{
\begin{array}{rcll}
    \Delta u_{\DD} + k^2 u_{\DD} & = & 0 & \text{ in } \D \\
    \partial_{\nn} u_{\DD} - ik_2u_{\DD} & = & h & \text{ on } \Gamma_{\mathrm{c}} \\
    u_{\DD} & = & u_{\mathrm{obs}} & \text{ on } \Gamma_{\mathrm{m}} \\
    \left[u_{\DD}\right] \; = \; [\partial_{\nn} u_{\DD}] & = & 0 & \text{ on }  \partial \Od_1, \\
\end{array}
\right.  &  
\left\{
\begin{array}{rcll}
    \Delta u_{\RR} + k^2 u_{\RR} & = & 0 & \text{ in } \D \\
    \partial_{\nn} u_{\RR} - ik_2u_{\RR} & = & h & \text{ on } \partial \D \\
    \left[u_{\RR}\right] \; = \; [\partial_{\nn} u_{\RR}] & = & 0 & \text{ on } \partial \Od_1.
    \end{array}
\right.
\end{array}\end{equation} 
The cost functional $\mathcal{K}$ has two parameters $\Od_1$ and $k_1$ that correspond to our two unknowns. Notice that, by the proposition above, the equations defining $u_{\DD}$ and $u_{\RR}$ are well posed. \\

We see that the solutions to the inverse problem minimize the functional as they are zeros of $\mathcal{K}$. However, we do not know if all minimizers are solutions to the inverse problem. This question is related to the identifiability of the problem: \textit{is knowing~$h$ and~$g$ enough to determine $\Od_1$ and $k_1$?} In our case, if we consider the three dimensional case $d=3$, then identifiability holds if all observation and scattering directions are known. We write the proposition as follows.
First, for all $(\Od_1 ,k_1) \in \mathcal{U}_{\rm adm} \times \mathfrak{K}_{\rm adm}$, we introduce the \textit{Robin-to-Dirichlet} operator $\Lambda_{\Od_1,k_1}$ defined by \begin{equation} \label{def:RobToDirOp}
\fonction{\Lambda_{\Od_1,k_1}}{\HH^{-1/2}(\partial \mathcal{D})}{\HH^{1/2}(\mathcal{\partial \D})}{h}{\gamma(u_\RR),}
\end{equation}
where $\gamma$ denotes the trace function on $\partial \D$ and $u_{\RR} \in V$ is the solution of~$\eqref{eq:urud}$.
 
\begin{proposition}\label{prop:uniqueness}
    Let $d=3$.
    If two couples of data $(\Od_{\rm a}, k_{\rm a})$ and~$(\Od_{\rm b}, k_{\rm b})$ in~$ \mathfrak{U}_{\rm adm} \times  \mathfrak{K}_{\rm adm}$ have the same Robin-to-Dirichlet operator, then they are equal. In other words, we have $$\Lambda_{\Od_{\rm a}, k_{\rm a}} = \Lambda_{\Od_{\rm b}, k_{\rm b}} \Longrightarrow \Od_{\rm a} = \Od_{\rm b} \text{ and } k_{\rm a} = k_{\rm b}.$$
\end{proposition}

\begin{remark}
Up to our knowledge, such a result for $d = 2$ remains an open problem in the present context. Nevertheless we highlight that similar identifiability results have been obtained in all dimensions in other contexts. Firstly, let us mention the work of Bourgeois {\it et al.}~\cite{Bourgeois2012}, where a perfectly electric conductor with a thin-layer is considered. Secondly, in the context of bounded elliptic problems, we point out the property~\cite[property 4.4]{Darde2010} where a single measurement is needed. 
\end{remark}

In the following, we aim to use an iterative Nesterov scheme (see section 4 for more details) to minimize the previous functional. To do so, we need to compute a gradient at each iteration with respect to both the shape $\Od_1$ and the wavenumber $k_1$. To calculate the shape gradient, we use the classical notion of shape derivative, that we briefly recall here (see~\cite{henrotpierre2018}).

\begin{definition}[Shape derivative]
    Let $J: \mathcal{D}_{adm} \rightarrow \mathbb{R}$ be a functional defined on a class of admissibles shapes $\D_{adm}$. Let $\mathcal{W}$ be a functional subspace of $W^{1,\infty}(\bR^d)$ and denote the identity $\bR^d \rightarrow \bR^d$ by $I$. The shape derivative of $J$ at $\Od$ is the Fréchet derivative of $$\begin{cases}
    \mathcal{W} \rightarrow \mathbb{R} \\
    \theta \mapsto J((I + \theta)(\Od))
    \end{cases}$$ if it exists. It will be denoted in the following by $\mathrm{D}_{\Od_1}\mathcal{J}(\Od)(\theta).$ \\
    
    In this article, the functional subspace $W$ used in the proof is $W_0^{1,\infty}(\D)$. It is the subspace of $W^{1,\infty}(\bR^d)$, made of functions with traces that vanish on $\partial \D$. \end{definition}
\subsection{Deterministic case}

We begin by giving the expression of the derivatives of $\mathcal{K}$ with respect to both the wavenumber $k_1$ and the shape $\Od_1$ in the deterministic case. 

\begin{proposition}[Derivative with respect to $\Od_1$]\label{prop:adjointfirststrategy}
Let $k_1 \in \mathfrak{K}_{\rm adm}$. For every $\Od_1 \in \mathcal{U}_{\rm adm}$, and $\theta \in \WW_0^{1,\infty}(\D)$, the shape derivative of the Kohn-Vogelius functional can be expressed as follows 

    \[
    \DD_{\Od_1}\mathcal{K}(\Od_1,k_1)(\theta) = \mathfrak{Re}\left(\int_{\partial \Od_1} (k_1^2 - k_2^2) (u_{\DD} \rho_{\DD} - u_{\RR} \rho_{\RR}) (\theta \cdot \nn) \right),
    \]
where $\rho_{\DD} \in V$ and $\rho_{\RR} \in V$ are solutions to the adjoint equations
\begin{equation}\label{eq:adjoints}
\begin{array}{c c}
\left\{
\begin{array}{rcll}
    \Delta \rho_{\DD} + k^2 \rho_{\DD} &=& -\bar{k}^2 (\overline{u_{\DD} - u_{\RR}}) & \text{ in } \D \\
    \partial_{\nn} \rho_{\DD} - ik \rho_{\DD} &=& \partial_{\nn} (\overline{u_{\DD}-u_{\RR})} & \text{ on } \Gamma_{\mathrm{c}} \\   
    \rho_{\DD} &=& 0 & \text{ on } \Gamma_{\mathrm{m}} \\
    \left[\rho_{\DD}\right] \; = \; \left[ \partial_{\nn} \rho_{\DD} \right] & = & 0 & \text{ on }  \partial \Od_1,\\
\end{array}
\right. &
\left\{
\begin{array}{rcll}
    \Delta \rho_{\RR} + k^2 \rho_{\RR} &=& -\bar{k}^2 (\overline{u_{\DD} - u_{\RR}}) & \text{ in } \D\\
    \partial_{\nn} \rho_{\RR} - ik\rho_{\RR} &=& \partial_{\nn}(\overline{u_{\DD}-u_{\RR}}) & \text{ on } \partial \D \\
    \left[\rho_{\RR}\right] \; = \; \left[ \partial_{\nn} \rho_{\RR} \right] & = & 0 & \text{ on }  \partial \Od_1. \\
\end{array}
\right.
\end{array}\end{equation}

\end{proposition}

\begin{proposition}[Derivative with respect to $k$]\label{prop:derivativekfirststrategy}
Let $\Od_1 \in \mathcal{U}_{\rm adm}$. For every $k_1 \in \mathfrak{K}_{\rm adm}$, the derivative of the Kohn-Vogelius functional \eqref{def:K} with respect to $k_1$ exists and

\[
    \frac{\partial \mathcal{K}(\Od_1,k_1)}{\partial k_1} = 2\int_{\Od_1} (\rho_D u_D - \rho_R u_R) k_1,
\]

where $u_{\DD} \in V$ and $u_{\RR} \in V$ are defined in the equations \eqref{eq:urud} and $\rho_{\DD} \in V$ and $\rho_{\RR} \in V$ are solutions to the previous equations \eqref{eq:adjoints}.

\end{proposition}

\subsection{Case with uncertainties}

The method detailed above is considered under ideal conditions, i.e., a deterministic case where the quantities $u_{\mathrm{obs}}$ and $h$ are known with precision. In practice, however, this is not the case. Two factors compromise the effectiveness of our approach. First, there is a noise on $u_{\mathrm{obs}}$ which arise from measurements. Second, $h$ depends on the incident wave vector that may deviate from our exact expectations. \\

 In order to take into account those uncertainties, we consider a stochastic framework where we assume that the errors follow a statistical law. Precisely, we let $(\Omega, \mathcal{T}, \bP)$ be a probability space, $h \in \LL_{\bP}^2(\Omega,\HH^{-1/2}(\partial \D))$ and $u_{\mathrm{obs}} \in \LL_{\bP}^2(\Omega,\HH^{1/2}(\Gamma_{\rm m}))$. For a given $\omega \in \Omega$, the inverse problem in the first strategy can be rewritten as finding $\Od_1 \in \mathcal{U}_{\rm adm}$, $k_1 \in \mathcal{K}_{\rm adm}$ and $u(\omega) \in V$ such that

\[
\left\{
\begin{array}{lll}
    \Delta u(\omega) + k^2 u(\omega) = 0 & \text{ in } \D \\
    \partial_{\nn} u(\omega) - iku(\omega) = h(\omega) & \text{ on } \partial \D \\
    u(\omega) = u_{\mathrm{obs}}(\omega) & \text{ on } \Gamma_{\rm m} \\
    \left[ u(\omega) \right] = \left[ \partial_{\nn} u(\omega) \right] = 0 & \text{ on }  \partial \Od_1.
\end{array}
\right. 
\]

 In practice, the main difficulty lies in the fact that we know only the realizations of the random variables~$u_{\mathrm{obs}}$ and~$h$. There are many ways to solve the problem with uncertainties. We could use a Monte-Carlo method or a stochastic gradient descent, but it would lead to very long computation time. To overcome this problem, we choose to minimize a combination of the expected value and the variance of the considered functional~$\mathcal{K}$, that is to say we are considering $$\tilde{\mathcal{K}}(\Od_1,k_1,\omega) = (1-\alpha) \bE(\mathcal{K}(\Od_1,k_1,\omega)) + \alpha \mathbb{V}(\mathcal{K}(\Od_1,k_1,\omega)),$$ for $\alpha \in [0,1]$. Hence, the problem is reduced to a purely deterministic problem which can be treated as before. For this purpose, we rely on a truncated Karhunen-Loève expansion (see~\cite{Faverjon2017,Dambrine2022} for more details). Let $M > 0$ be an integer, $Y$ be a centered and normalized random variable with all the moments of order lower than $4$ finite. We consider $M$ independent and identically distributed random variables $Y_i \sim Y$, $(g_m)_{1 \le m \le M} \in\HH^{1/2}(\Gamma_{\rm m})$ and $(h_m)_{1 \le m \le M} \in\HH^{-1/2}(\Gamma_{\rm m})$ such that
\begin{equation}\label{expKL}
    u_{\mathrm{obs}}(x,\omega) = g_0(x) + \sum_{m=1}^M g_m(x)Y_m(\omega) \quad \text{and} \quad h(x,\omega) = h_0(x) + \sum_{m=1}^M h_m(x)Y_m(\omega).
\end{equation}

By linearity of the Helmholtz equation, we get

\[u_{\RR}(\omega) = u_{\RR_0} + \sum_{m = 1}^M u_{\RR_{\rm m}} Y_{\rm m}(\omega) \quad \text{and} \quad u_{\DD}(\omega) = u_{\DD_0} + \sum_{m = 1}^M u_{\DD_{\rm m}} Y_{\rm m}(\omega),\]

 where $u_{\DD_m} \in V$ and $u_{\RR_m} \in V$ solve for all $0 \le m \le M$,

\[
\begin{array}{c c}
\left\{
\begin{array}{lll}
    \Delta u_{\DD_m} + k^2 u_{\DD_m} = 0 & \text{ in } \D \\
    \partial_{\nn} u_{\DD_m} - iku_{\DD_m} = h_m & \text{ on } \Gamma_{\mathrm{c}} \\
    u_{\DD_m} = g_m & \text{ on } \Gamma_{\mathrm{m}} \\
    \left[u_{\DD_m}\right] = \left[ \partial_{\nn} u_{\DD_m} \right] = 0 & \text{ on }  \partial \Od_1, \\
\end{array}
\right. &
\left\{
\begin{array}{lll}
    \Delta u_{\RR_m} + k^2 u_{\RR_m} = 0 & \text{ in } \D\\
    \partial_{\nn} u_{\RR_m} - iku_{\RR_m} = h_m & \text{ on } \partial \D \\
    \left[u_{\RR_m}\right] = \left[ \partial_{\nn} u_{\RR_m} \right] = 0 & \text{ on }  \partial \Od_1. \\
\end{array}
\right.
\end{array}\] We can compute $\mathbb{E}(\mathcal{K}(\Od_1,k_1,\omega))$ and $\mathbb{V}(\mathcal{K}(\Od_1,k_1,\omega))$ in any cases, but we will consider the case $Y \sim \mathcal{N}(0,1)$ to simplify the expressions of the variance in the following. 

\begin{proposition}[Expression of the expected value of the functional]\label{prop:expectedvalue}

In the framework described above, we have
\[\mathcal{E}_{\mathcal{K}}(\Od_1,k_1) = \bE(\mathcal{K}(\Od_1,k_1,\omega)) = \frac{1}{2} \int_{\D} \sum_{m= 0}^M |\nabla(u_{\DD_m} - u_{\RR_m})|^2.\]

For $\theta \in \mathrm{W}^{1,\infty}_0(\bR^d)$, $\Od_1 \in \mathcal{U}_{\rm adm}$ and $k_1 \in \mathfrak{K}_{\rm adm}$, the derivatives given by the adjoint exist and are given by

\[\DD_{\Od_1}\mathcal{E}_{\mathcal{K}}(\Od_1,k_1)(\theta) = \sum_{m = 0}^M \mathfrak{Re}\left(\int_{\partial \Od_1} (k_1^2 - k_2^2)( u_{\DD_m} \rho_{\DD_m} - u_{\RR_m} \rho_{\RR_m})(\theta \cdot \nn) \right),\]

\begin{equation}\label{Kprime}
\frac{\partial \mathcal{E}_{\mathcal{K}}(\Od_1,k_1)}{\partial k_1} = 2 \sum_{m=0}^M \int_{\Od_1} (\rho_{\DD_{\mathrm{m}}} u_{\DD_{\mathrm{m}}} - \rho_{\RR_{\mathrm{m}}}u_{\RR_{\mathrm{m}}})k_1,
\end{equation}
where $\rho_{D_m} \in V$ and $\rho_{R_m} \in V$ solve, for all $1 \le m \le M$,

\[
\begin{array}{c c}
\left\{
\begin{array}{lll}
    \Delta \rho_{\DD_m} + k^2 \rho_{\DD_m} = -\overline{k}^2 (\overline{u_{\DD_m} - u_{\RR_m}}) & \text{ in } \D \\
    \partial_{\nn} \rho_{D_m} - ik \rho_{\DD_m} = \partial_{\nn} (\overline{u_{\DD_m}-u_{\RR_m})} & \text{ on } \Gamma_{\mathrm{c}} \\   
    \rho_{\DD_m} = 0 & \text{ on } \Gamma_{\mathrm{m}} \\
    \left[\rho_{\DD_m}\right] = \left[\partial_{\nn} \rho_{\DD_m}\right] = 0 & \text{ on }  \partial \Od_1, \\
\end{array}
\right. &
\left\{
\begin{array}{lll}
    \Delta \rho_{\RR_m} + k^2 \rho_{\RR_m} = -\overline{k}^2 (\overline{u_{\DD_m} - u_{\RR_m}}) & \text{ in } \D \\
    \partial_{\nn} \rho_{\RR_m} - ik\rho_{\RR_m} = \partial_{\nn}(\overline{u_{\DD_m}-u_{\RR_m}}) & \text{ on } \partial \D \\
    \left[\rho_{\RR_m}\right] = [\partial_{\nn} \rho_{\RR_m}] = 0 & \text{ on }  \partial \Od_1.
\end{array}
\right.
\end{array}\]

\end{proposition}
\begin{proposition}[Expression of the variance of the functional]\label{prop:variance}
We introduce the notation $\mathrm{d}_i = u_{\RR_i} - u_{\DD_i}$ and $\mathrm{p}_{j,i} = u_{\DD_i}\rho_{\DD_j} - u_{\RR_i}\rho_{\RR_j}$. Let $\Od_1 \in \mathcal{U}_{\rm adm}$, $k_1 \in \mathfrak{K}_{\rm adm}$. In the framework described above, if $Y \sim \mathcal{N}(0,1)$, we have
\begin{multline}
    \mathcal{V}_{\mathcal{K}}(\Od_1,k_1) = \mathbb{V} \left(\mathcal{K}(\Od_1,k_1,\omega) \right) = \frac{1}{4} \left(\int_{\D} |\nabla \mathrm{d}_0|^2\right)^2 + \frac{1}{2} \left(\int_{\D} |\nabla \mathrm{d}_0|^2\right)  \left( \sum_{m = 1}^M \int_{\D}|\nabla \mathrm{d}_m|^2\right)
    \\ + \sum_{m = 1}^M \left( \int_{\D} \mathfrak{Re}\left( \nabla \mathrm{d}_0 \cdot \nabla \overline{\mathrm{d}_m}\right)\right)^2 +  \mathfrak{Re}\left( \sum_{1 \le m < l \le M} \left(\int_{\D} \nabla \mathrm{d}_m \overline{\nabla \mathrm{d}_l}\right)^2 \right) 
    + \frac{1}{2} \sum_{m = 1}^M \left(\int_{\D} |\nabla \mathrm{d}_m|^2\right)^2.
\end{multline}

 For $\theta \in \WW_0^{1,\infty}(\D)$, the derivatives given by the adjoint method are now \[\mathrm{D}_{\Od_1} \mathcal{V}_{\mathcal{K}} \left((\Od_1,k_1) \right)(\theta) = \mathfrak{Re}\left( \int_{\partial \Od_1} (k_1^2-k_2^2) \mathcal{G} (\theta \cdot n)\right) \quad \mathrm{and} \quad \frac{\partial \mathcal{V}_{\mathcal{K}}(\Od_1,k_1)}{\partial k_1} = 2 \int_{\Od_1} \mathcal{G} k_1\]where 
\begin{multline}
    \mathcal{G} = \frac{1}{2} \left(\int_{\D} |\nabla \mathrm{d}_0|^2\right)\sum_{m = 1}^M \mathrm{p}_{m,m} + 2\sum_{m = 1}^M \left( \int_{\D} \mathfrak{Re}\left( \nabla \mathrm{d}_0 \cdot \nabla \overline{\mathrm{d}_m}\right)\right)\mathrm{p}_{0,m} \\
    + 4 \sum_{1 \le m < l \le M} \mathfrak{Re}\left(\int_{\D} \nabla \mathrm{d}_m \overline{\nabla \mathrm{d}_l}\right)  \mathrm{p}_{m,l} + \sum_{m = 1}^M \left(\int_{\D} |\nabla \mathrm{d}_m|^2 \right) \mathrm{p}_{m,m}. 
\end{multline}

\end{proposition}

\section{Numerical experiments}
\label{sec:numexp}
To illustrate the effectiveness of our approach, we present numerical experiments in 2D. 

\subsection{Algorithm}
To ensure visibility of the unknown target, we consider different incidents directions $\theta_j$ for $1 \le j \le J$. For each direction, we take $$h^j= i(k - k_{\rm inc}^j(\theta_j) \cdot n) u_{i_0} \mathrm{exp}(-ik_{\rm inc}^j(\theta_j)\cdot r),$$ where $k_{\rm inc}^j$ corresponds to the incident wave vector given by $k_{\rm inc}^j = -k(\cos(\theta_j),\sin(\theta_j))$. It is important to note that we have in total $J$ functionals rather than a single one. We chose to minimize each functional successively over a given number of iterations. Our method consists of a sequential process: we begin by minimizing with respect to the first direction of the plane wave, then the second, and so on, until all directions have been considered. This cycle is then repeated, starting again with the first direction. Notice that Furuya and Potthast investigate in \cite{Furuya2022} the minimization of the functional in such a sequential process for inverse scattering problem. \newline

The descent method used in the following is the Nesterov's accelerated gradient scheme~\cite{Nesterov1983}: 

\begin{equation*}
    \left\{
\begin{aligned}
  x_{n+1} &= y_n - \tau_{\Od_1,k_1} \nabla f(y_n) \\
  y_{n+1} &= x_{n+1} + \frac{\eta_t - 1}{\eta_{t+1}}\left(x_{n+1} - x_n\right), \quad y_0 = x_0 \\
  \eta_{t+1} &= \frac{1 + \sqrt{1 + 4\eta_t^2}}{2}, \quad \eta_0 = 1 ,
\end{aligned}\right.
\end{equation*} with restarting~\cite{Donoghue2015}: if $\mathcal{K}(x_{n+1}) > \mathcal{K}(x_n)$, we take $\eta_0 = 1$, $x_0 = x_n$, $y_0 = x_n$. The variables in the algorithm (e.g. $x_n$ and $y_n$) are vectors of size $2$ of the form $[\Od_1,k_1]$ and the step time $\tau$ is $\tau_{\Od_1,k_1}=[\tau_{\Od_1},\tau_{k_1}]$.  Hence, both shape and parameters are updated at each iteration at the same time, with their own gradient and step $\tau$. \\

The choice of the two descents step is a key point of our algorithm. If the step size are too large, the process becomes unstable, potentially leading to divergence or even failure to converge at all. On the other hand, if the step sizes are too small, the reconstruction becomes excessively slow. \\

In addition, we consider a truncated Fourier series parametrization for the boundary of the object $\mathcal{O}_1$:

\[\partial \Od_1 = \begin{bmatrix}
    x_c \\
    y_c
\end{bmatrix} +\left(r + \sum_{n = 1}^N a_n \mathrm{cos}(n\theta) + b_n\mathrm{sin}(n\theta) \right) \begin{bmatrix}
     \mathrm{cos}(\theta)\\
     \mathrm{sin}(\theta)
\end{bmatrix}.\]

To avoid the instability of the inverse problem, the Fourier coefficients are introduced incrementally. The first five iterations are used to adjust the center. The radius is introduced after this fifth iteration, and each subsequent Fourier coefficient is added at multiples of $30$ iterations. Hence, the more iterations we do, the bigger is the number of Fourier coefficients we accumulate by the end of our algorithm. Note that we could also choose to add them following a criterion, such as adding some coefficients if the functional is not changing over a certain number of iterations. Nonetheless, the goal here is to validate the method. The general idea of the proposed scheme is summed up in the pseudo-code~\ref{alg:desNest}. 

\begin{algorithm}
\caption{Descent algorithm}\label{alg:desNest}
\begin{algorithmic}
\Require $N_{\mathrm{iter}}, N_{\mathrm{seq}}$, $\theta_{\rm inc}$
\For{$0 \le n < N_{\mathrm{iter}}$}
\For{$\theta_j \in \theta_{\rm inc}$}
\State Nesterov scheme with $N_{\mathrm{seq}}$ iterations and the functional $\mathcal{K}^j$:
\State update of the wavenumber and the Fourier coefficients
\If{n/30 is an integer}
\State add one Fourier mode
\EndIf
\EndFor
\EndFor
\end{algorithmic}
\end{algorithm}

\begin{remark}
Note that a new mesh is created at each iteration.
\end{remark}

\subsection{Implementation}

We conduct numerical experiments to assess the effectiveness of the different methods developed earlier in section~\ref{sec:mainres}. We use the PDE solver FreeFEM++ (see~\cite{Hecht2012}) and the package pyfreefem (see~\cite{Feppon2020}) to interface with Python. The FreeFEM++ version we use is parallelized using the MPI library. The resolution of the inverse problem requires a high-performance computing (HPC) cluster, specifically Pyrene developed at the Université de Pau. The forward problem is first modeled using P2 finite elements on a highly refined mesh. The solution is then computed on a new mesh with P1 finite elements to avoid the inverse crime problem. \newline

We use the parallelization for the robust case only. It is implemented as a process distribution, where each process corresponds to a term in the Karhunen–Loève expansion~\eqref{expKL}. Thanks to this method, we were able to increase the order of the expansion while maintaining approximately the same computation time. With~$5$~terms in the expansion, each iteration requires approximately $60$s.

\subsection{Results} 

\begin{remark}\label{remark:zone}
For a given object of characteristic size $r_{\rm obj}$, compared to the wavelength $\lambda = c_0/(f\mathfrak{Re}(\sqrt{\varepsilon}))$ (where $c_0$ denotes the celerity in the vacuum, $f$ the frequency and $\varepsilon$ the permittivity of the surrounding medium), the behavior of the scattered field can be decomposed into three different regimes. When $r_{\rm obj} \ll \lambda$ (typically~$r_{\rm obj} \lesssim \lambda/10$), corresponding to the so-called {\it Rayleigh region}, the scattered field depends mainly on the global size of the object and in a lesser amount to its average permittivity contrast, with very weak sensitivity to its detailed shape, making shape reconstruction severely ill-posed. When $r_{\rm obj}$ becomes comparable to $\lambda$, corresponding to the so-called {\it Mie region}, the scattered field exhibits a strong oscillatory behavior with both the object geometry and the permittivity contrast: this regime thus contains more information, but the oscillations generally induce multiple local minima. Finally, when $r_{\rm obj} \gg \lambda$ (in general~$r_{\rm obj} \gtrsim 10\lambda$), corresponding to the so-called {\it Fresnel region}, scattering is predominantly controlled by boundary geometry, i.e., the shape, reflection, and absorption effects, i.e., mostly the imaginary part of $\varepsilon$: this also highlight the need of dual methods that can recover the information in different part of this spectrum.
\end{remark}

In all the experiences below, we consider a general domain $\mathcal{D}$, which is modeled as a square with length~$2$. We consider four incident waves coming from directions $\theta_i = 0, \frac{\pi}{2}, \pi$ and $\frac{3\pi}{2}$, with each having the same amplitude of $1$ and the same frequency $f=30~$MHz. 
The space of measures $\Gamma_{\rm m}$ is constituted of four intervals of length $2/3$ centered on each side of the boundary of $\D$. 
In the following, our main objective is to reconstruct a shape that does not enter into our Fourier parametrizations. We will name it the "kite" and it can be parameterized by  

$$
\begin{bmatrix}
x \\
y
\end{bmatrix} = \begin{bmatrix}
- 0.26 + 0.4 \cos(t) + 0.52 \cos^2(t) \\
0.6 \sin(t)
\end{bmatrix} , \qquad t \in [0,2\pi).
$$ 
For the object we consider a wavenumber, $k_1$, associated to a permittivity $\varepsilon_1 = 50 + 29i$. The outside permittivity is taken $\varepsilon_2 = 10 + 4i$. For all our tests,  the initial guess for the shape is a disk with a radius of~$0.6$.\\

\begin{remark}
Considering the characteristic size of the object, $r_{\rm obj} \simeq \lambda/6$, 
this means we are in the Mie region. We thus expect our proposed framework to provide accurate results while significant oscillations may occur in the gradient descent. In addition, since the contrast is high here, the oscillations typical of the Mie region would be amplified, further increasing the difficulty of the underlying inverse scattering problem. The relatively large imaginary part of the permittivity also leads to stronger attenuation of the scattered field, which reduces the available information and further complicates the reconstruction in all the tests considered.\end{remark}

\subsubsection{Deterministic reconstruction}
\paragraph{Deterministic reconstruction of the shape at fixed wavenumber.} 
 
\begin{figure}[htbp]
  \centering
    \subfloat[Evolution of the functional]{\includegraphics[trim={1cm 0 2cm 1cm}, clip, width=0.51\textwidth]{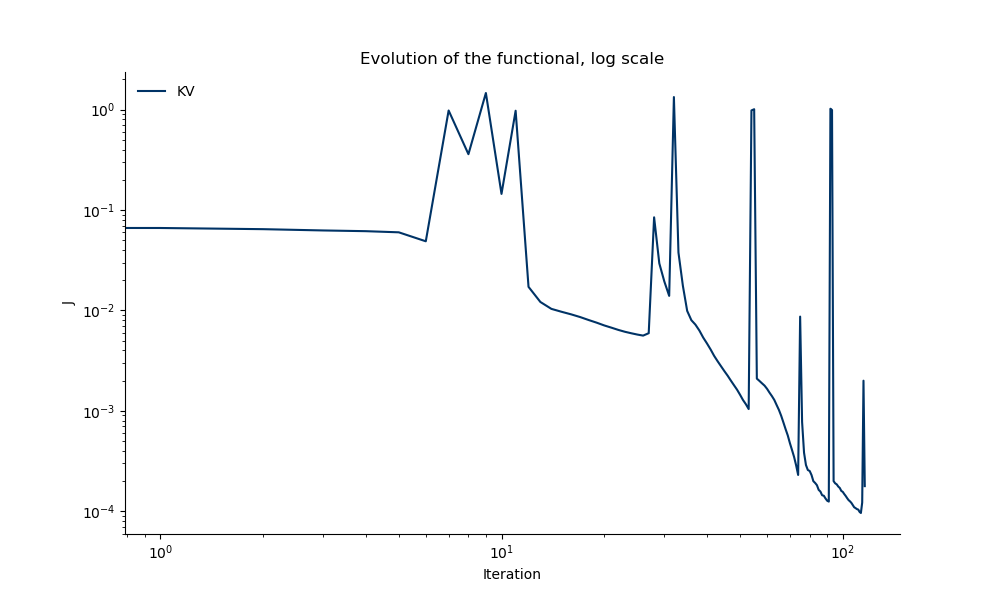}}
  \hfill
  \subfloat[Retrieved shape]{
    \includegraphics[trim={0 1.5cm 0 1.5cm},width=0.48\textwidth, clip]{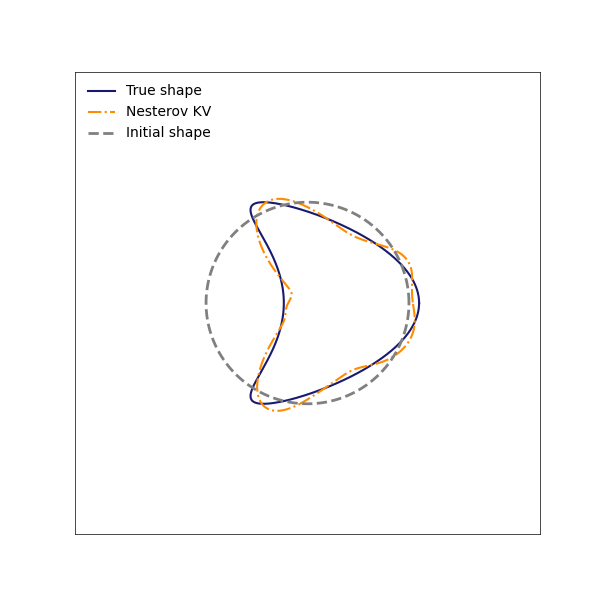}}
  \caption{Reconstruction of the shape at fixed wavenumber}
  \label{fig:deterministshape}
\end{figure}

First, we aim to reconstruct the shape assuming that we know the wavenumber. Hence, we take $k_1$ as the initial wavenumber and attempt to reconstruct the kite shape. The results, depicted in Figure~\ref{fig:deterministshape}, show the final retrieved shape (orange, dashed) closely approximately the true kite geometry (plain). In Figure~\ref{fig:deterministshape}, we also represent the evolution of the functional. The initial shape is represented in dotted. All the figures represented will have the same convention. The descent is performed with a total of $400~$iterations, with $N_{\mathrm{iter}} = 4$ descent directions, $N_{\mathrm{seq}} = 25$ sequential step per direction, and a step size of \textbf{$0.02$}. \\

The objective functional does not decrease monotonically. Indeed, the functional we minimize along the descent depends on the incident direction. Hence, we can observe some jumps in the functional decrease. They correspond to the changes of initial direction. Moreover, the difference between the total number of iterations (400) and the number of points plotted is due to the removal of certain iterations at restart points within the algorithm. Nevertheless, we can see that the functional decreases globally and converges to $0$ with the number of iterations. All the other experiments have a similar behaviour for the functional. 

\begin{figure}[h!tbp]
    \centering
    \subfloat[Initialization with $0.9k_1$]{\includegraphics[trim={1.5cm 0 2.5cm 0cm}, clip, width=0.48\textwidth]{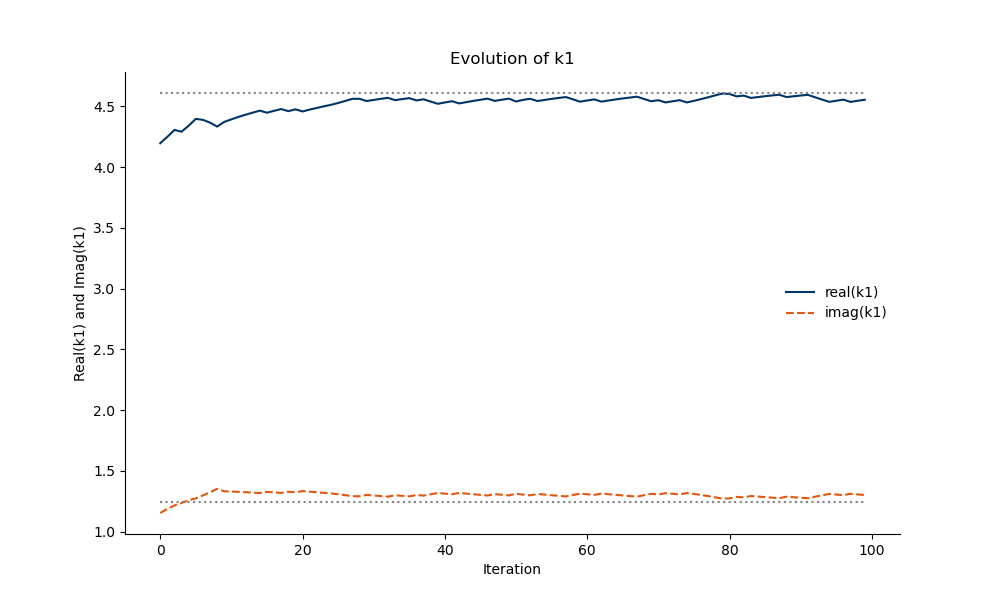}}
    \hfill
    \subfloat[Initialization with $1.2k_1$]{
        \includegraphics[width=0.48\textwidth, trim={1.5cm 0 2.5cm 0cm}, clip]{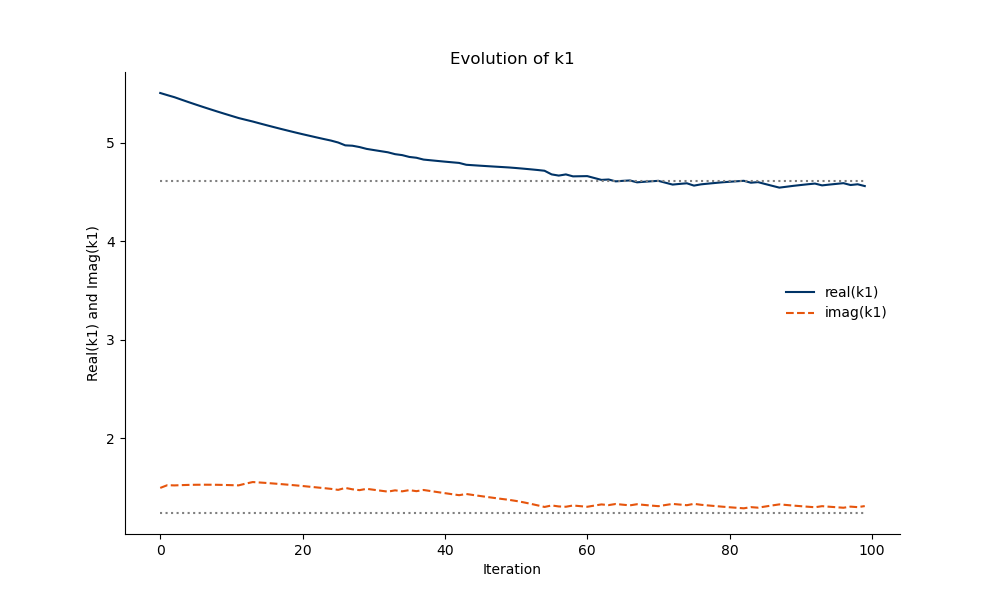}}
    \caption{Reconstruction of the wavenumber at fixed shape}
    \label{fig:deterministparameter}
\end{figure}

\paragraph{Deterministic reconstruction of the wavenumber at fixed shape.} We also attempted to reconstruct the wavenumber $k_1$ using the unknown as the initial shape. Initially, we set the shape update step $\tau_{\Od_1}$ to $0$, so that the algorithm is a Nesterov descent on the parameter $k$ only. In this experience only, the objective shape is the initial shape: the disk of radius $0.6$. Also, we have chosen $N_{\rm iter} = 1$ and $N_{\mathrm{seq}} = 25$. However, this intuitive approach failed. We observed that the functional $J$ is extremely sensitive to the unknown shape (we refer to~\cite[Section~4.3]{Darde2013} where this observation was done in another context). Indeed retrieving the permittivity is a very ill-posed inverse problem in general. In particular, in our setting, i.e., the size of the object of around~$\lambda/6$, we are in the Mie diffraction zone where the scattering waves variations are highly oscillating with the area of the target and highly dependent on the difference of permittivities (see Remark~\ref{remark:zone}). Therefore, this sensitivity suggests that even very small shape variations of the shape (for instance, those introduced by copying the shape to a new mesh) result in large variations in $J$. Our solution was to consider a very small step for the shape update ($0.0001$) compared to the wavenumber update~($0.4$). This strategy keeps the shape very close to the true one. We tested this with two initial wavenumbers: $0.9k_1$ and $1.2k_1$. All the other parameters remain the same. In Figure~\ref{fig:deterministparameter}, we plot the evolution of the real and imaginary parts of the wavenumber during the reconstruction. \\

We can notice that the results are very accurate in both cases, even if no regularization is performed. Furthermore, this highlight the importance of the proposed dual strategy, i.e., incorporating the shape reconstruction, for the inverse scattering problem. It shall also be noted that obstacle detection and reconstruction at low frequency is of interest for the electromagnetic community~\cite{wen2013algorithm}.

\begin{figure}[h!tb]
  \centering
  \subfloat[$3~$MHz]{\includegraphics[trim={0 0.8cm 0 0}, width=0.48\textwidth, clip]{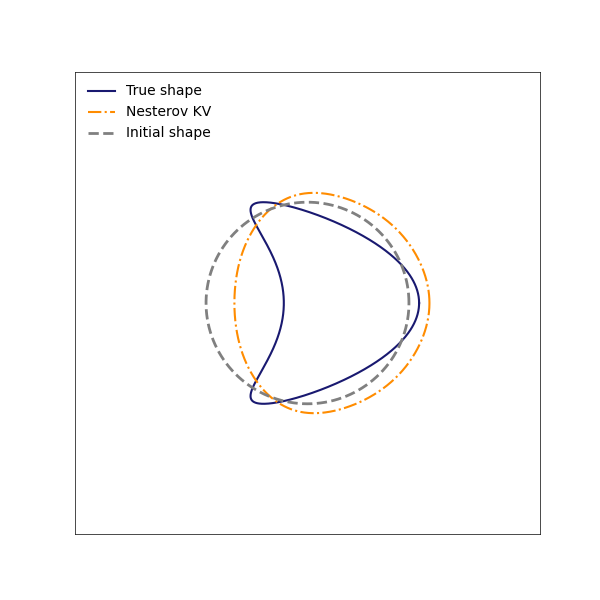}\label{fig:lowfrequency}} \hfill 
  \subfloat[$30~$MHz]{\includegraphics[trim={0 0.8cm 0 0}, width=0.48\textwidth, clip]{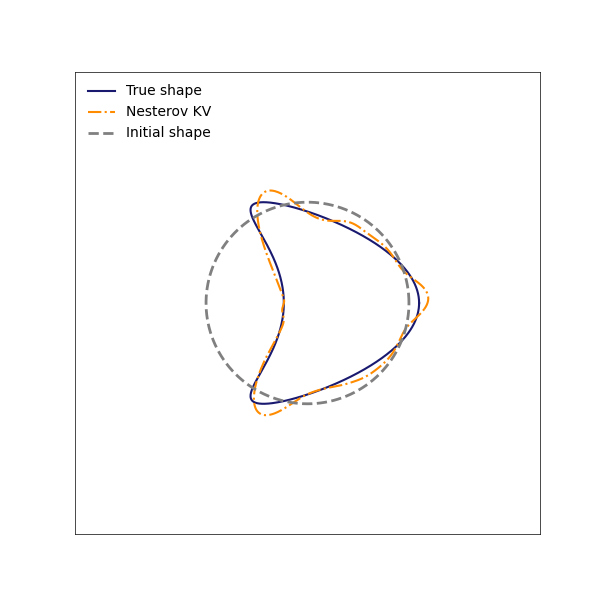}} \caption{Simultaneous reconstruction for two frequencies}
   \label{fig:determinist}
\end{figure}

\paragraph{Deterministic reconstruction of both wavenumber and shape.}
We now present a deterministic reconstruction when the two parameters $\Od_1$ and $k_1$ are unknowns. In Figure \ref{fig:determinist}, we represent the retrieved shape for two different frequencies. The initial values are $k = 0.6k_1$ and the disk of radius $0.6$. We performed the algorithm with $N_{\rm iter} = 3$ and $N_{\rm seq} = 25$ and the steps are $\tau_{\Od_1} = 3$ and $\tau_{k_1} = 4$ at $3~$MHz, and $\tau_{\Od_1} = 0.02$ and $\tau_{k_1} = 0.4$ at $30~$MHz. Even if the problem is very ill posed in the first case, we still expect our method to give good results. \newline

Firstly, notice that the reconstruction at $30$~MHz, i.e. in the Mie region, is efficient. Secondly,
with a frequency of $3$ MHz, we are in the Rayleigh diffraction zone (the size of the object is around~$\lambda/60$). Since, the functional is very sensitive to every modifications of the effective area, we have to put a step that makes a very small modification of the shape at each iteration. As a result, convergence is very slow. We observe, however, that the Fourier coefficients exhibit a small convexity aligning with the objective shape. Further investigation into the degrees of freedom of these coefficients may offer a path toward accelerating the process.

\subsubsection{Reconstruction with uncertainties}

\paragraph{Simultaneous reconstruction with uncertainties and the expected value.} We now attempt the simultaneous reconstruction of both the shape and the parameter with noisy data. To do so, $h_0$ is assumed to be the deterministic data, with $g_0$ being the associated observation $u_{\rm obs}$. Then, we consider $g_{m}$ and $h_m$ of the following forms: $$g_m = \beta_g \sin(m \phi) \quad \text{and} \quad h_m = \beta_h \sin(m \phi) i(k - k_{\rm inc}^j(\theta_j) \cdot n) \mathrm{exp}(-ik_{\rm inc}^j(\theta_j)\cdot r),$$ where $(\phi,r)$ are the usual polar coordinates. \newline

\begin{figure}[h!tb]
  \centering
  \subfloat[Reconstruction of the shape with $\beta_g = 5\%$ and $\beta_h = 1\%$]{
    \includegraphics[trim={0 0.5cm 0 0}, clip, width=0.48\textwidth]{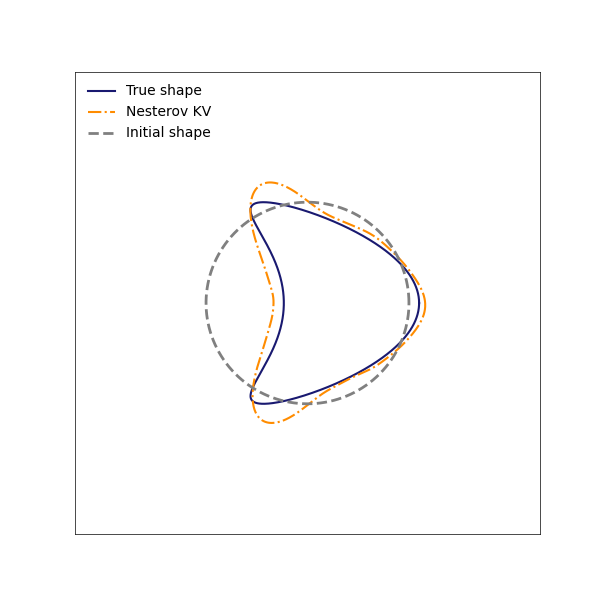}
  }
  \hfill
  \subfloat[Reconstruction of the shape with $\beta_g = 1\%$ and $\beta_h = 5\%$]{
    \includegraphics[trim={0 0.5cm 0 0}, clip, width=0.48\textwidth]{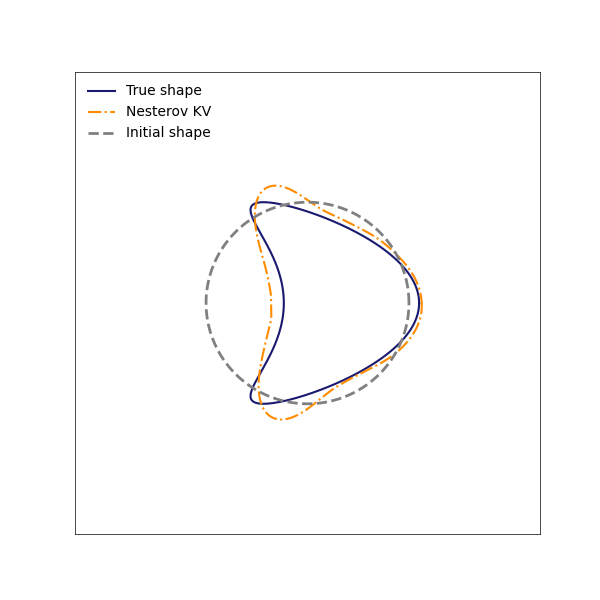}
  }

  \vspace{0.5em}

  \subfloat[Evolution of $k_1$ with $\beta_g = 5\%$ and $\beta_h = 1\%$]{
    \includegraphics[trim={1.5cm 0 3cm 0cm}, clip, width=0.48\textwidth]{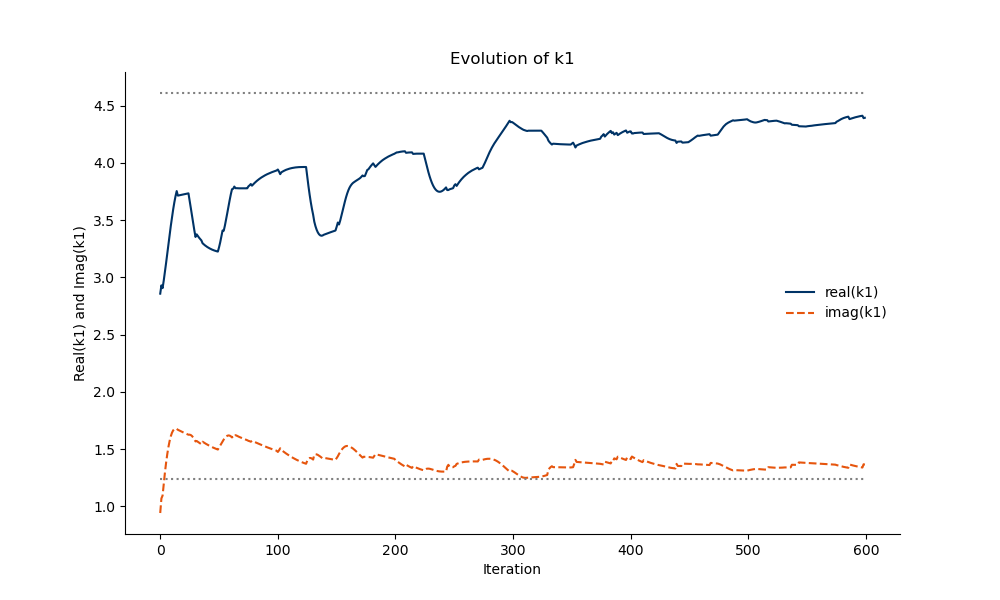}
  }
  \hfill
  \subfloat[Evolution of $k_1$ with $\beta_g = 1\%$ and $\beta_h = 5\%$]{
    \includegraphics[trim={1.5cm 0 3cm 0cm}, clip, width=0.48\textwidth]{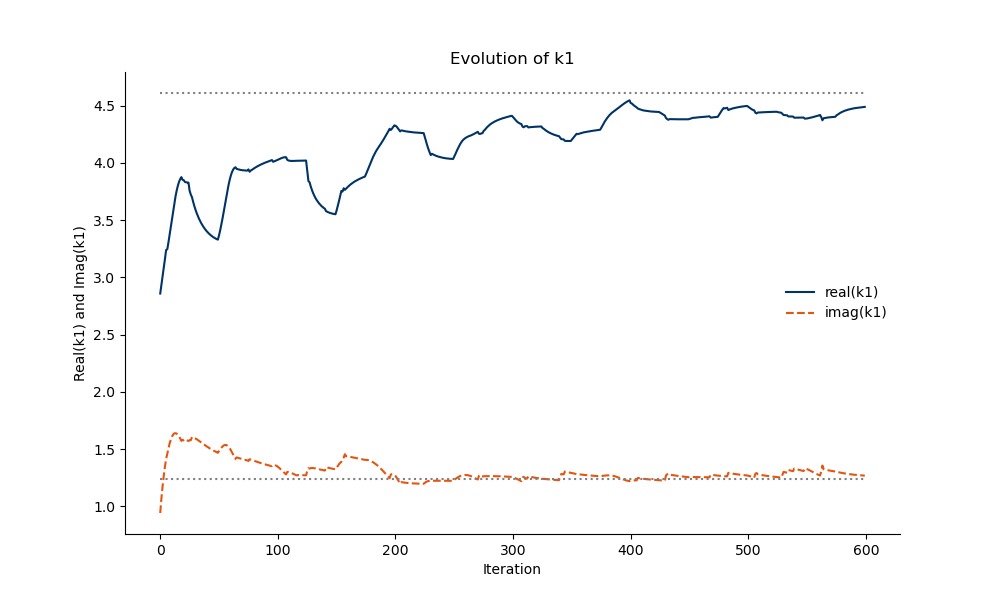}
  }

  \caption{Simultaneous reconstruction with noisy data and minimization of the expected value.}
  \label{fig:simultaneous}
\end{figure}

\begin{figure}[h!tb]
  \centering
  \subfloat[$\alpha = 0$]{
    \includegraphics[trim={0.75cm 1.25cm 0.75cm 1.5cm}, clip, width=0.32\textwidth]{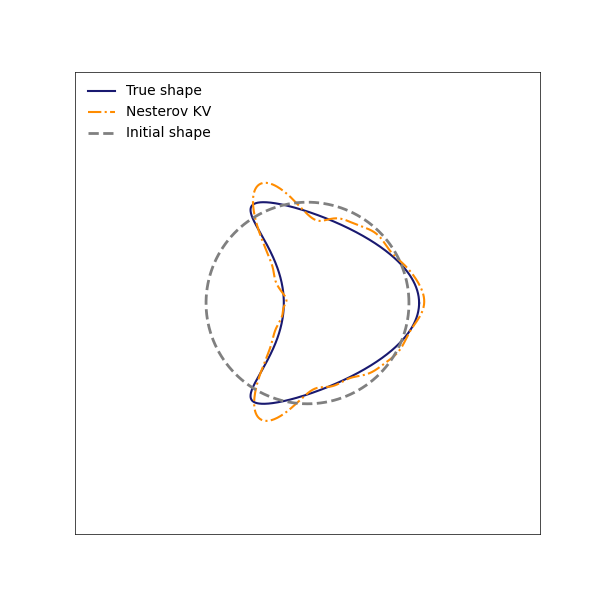}
  }
  \subfloat[$\alpha = 0.2$]{
    \includegraphics[trim={0.75cm 1.25cm 0.75cm 1.5cm}, clip, width=0.32\textwidth]{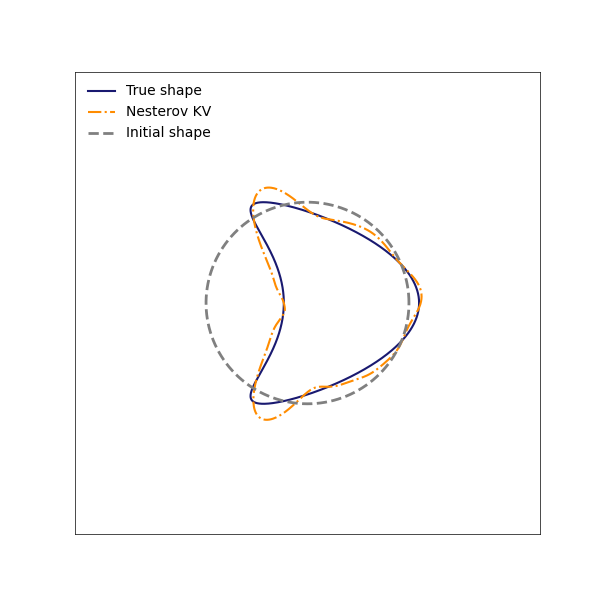}
  }
  \subfloat[$\alpha = 0.4$]{
    \includegraphics[trim={0.75cm 1.25cm 0.75cm 1.5cm}, clip, width=0.32\textwidth]{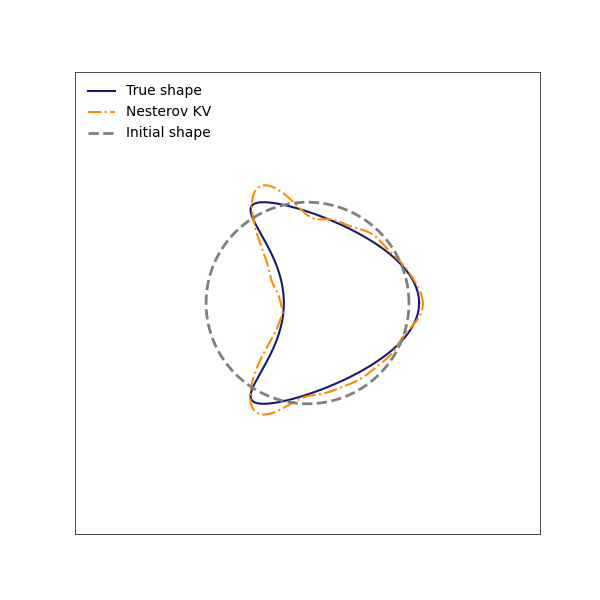}
  } \newline
  \subfloat[Evolution of $\mathfrak{Re}(k_1)$]{\includegraphics[scale=0.48, trim={1.5cm 0 0cm 0cm}, clip]{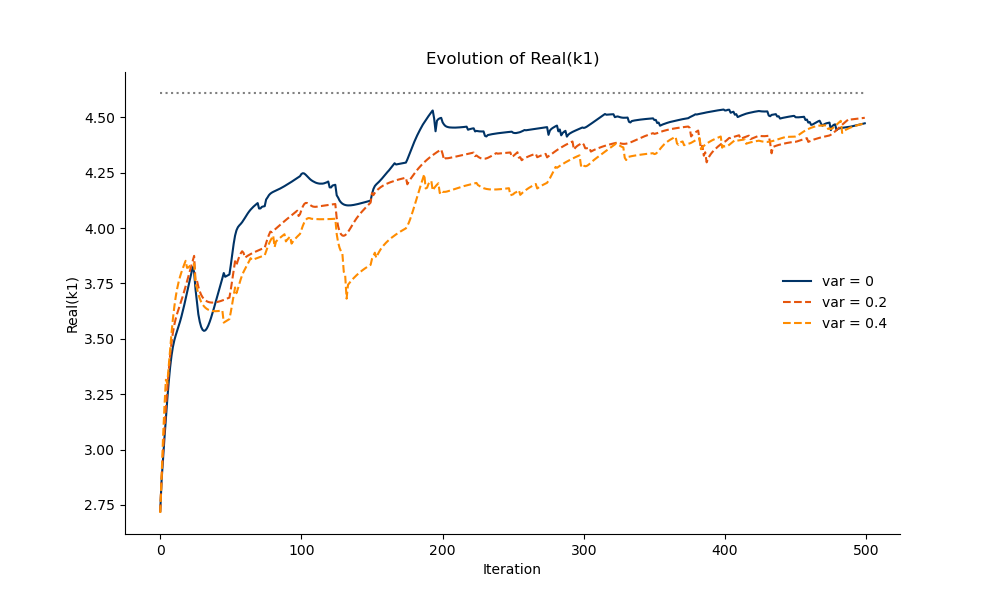}    \label{fig:robustvar}
}
  \caption{Simultaneous reconstruction with noisy data and a convex combination of the expected value and the variance}
  \label{fig:robustexpected}
\end{figure}

We first consider the minimization of the expected value only. We investigate one case with $5\%$ of noise on $g$ and $1\%$ on $h$, and conversely $1\%$ of noise on $g$ and $5\%$ on $h$, see \ref{fig:robustexpected}. By "percentage of noise," we refer to the noise level relative to the squared of the $\mathrm{L}^2$-norm on the boundary. Thus, $\beta_g$ and $\beta_h$ were chosen to match these expected noise percentages. As an example, $5\%$ of noise on $g$ and $1\%$ on $h$ correspond to
$$0.05^2 = \beta_g \bE\left (\frac{||g_1Y_1 + g_2Y_2 + g_3Y_3 + Y_4g_4||_{\mathrm{L}^2(\Gamma_{\rm m})}^2}{||g_0||_{\mathrm{L}^2(\partial \D)}^2}\right)$$
and
$$0.01^2 = \beta_h \bE\left(\frac{||h_1Y_1 + h_2Y_2 + h_3Y_3 + h_4Y_4||_{\mathrm{L}^2(\partial \D)}^2}{||h||_{\mathrm{L}^2(\partial \D)}^2}\right).$$ In Figure \ref{fig:simultaneous}, we represent the final shape and the evolution of the wavenumber for both experiments. 

With these levels of noise, we see that the method is able to retrieve both shape and wavenumber and remain robust under uncertainties. 

\paragraph{Simultaneous reconstruction with uncertainties by minimization of the expected value and variance.} Lastly, we put $3\%$ of errors on both $h$ and $g$ and consider $3$ different combinations of expected value and variance. Hence, the functional considered is $$\mathcal{J} = (1-\alpha)\mathcal{E}_{\mathcal{K}}(\Od_1,k_1) + \alpha \mathcal{V}_{\mathcal{K}}(\Od_1,k_1).$$ In Figure \ref{fig:robustexpected}, we represent the final shape for three values of $\alpha$. In Figure \ref{fig:robustvar}, we represent the evolution of the real parts of the wavenumber for the three experiments. \newline

We observe that increasing the value of $\alpha$ smooths the reconstructed shape, leading to a slightly worse reconstruction. On the other hand, the wavenumber becomes more stable with increasing variance. This tests show that the variance act as a regularizer for the robust approach and thus the choice of a good $\alpha$ is critical.

\subsection{Comparison with classical methods}

In this subsection, we aim to motivate the choices made above in this manuscript. In particular, we compare the Nesterov method with the gradient descent and the Kohn-Vogelius functional with the least squares functional.

\subsubsection{Gradient descent}
Nesterov showed in \cite{nesterov} that the method which now bears his name achieves a quadratic convergence rate~$O(1/N^2)$ for smooth convex functions, whereas the standard gradient descent only achieves a sublinear rate~$O(1/N)$. This scheme is known to be optimal among first-order methods. However, in the context of inverse problems that are highly non-convex and ill-posed, these theoretical convergence bounds do not strictly apply. \\

\begin{figure}[h!tb]
  \centering
  \subfloat[General comportment]{
    \includegraphics[trim={1cm 0 2cm 1cm}, clip, width=0.48\textwidth]{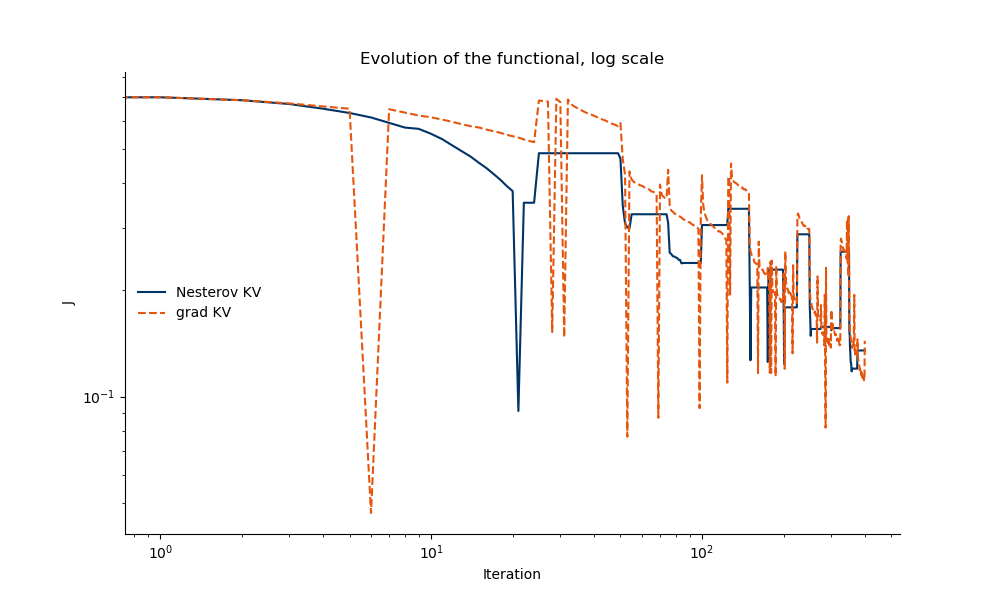}\label{fig:7a}
  } 
  \hfill
  \subfloat[Reconstruction of the wavenumber only]{
    \includegraphics[trim={1.5cm 0 2.5cm 0cm}, clip, width=0.48\textwidth]{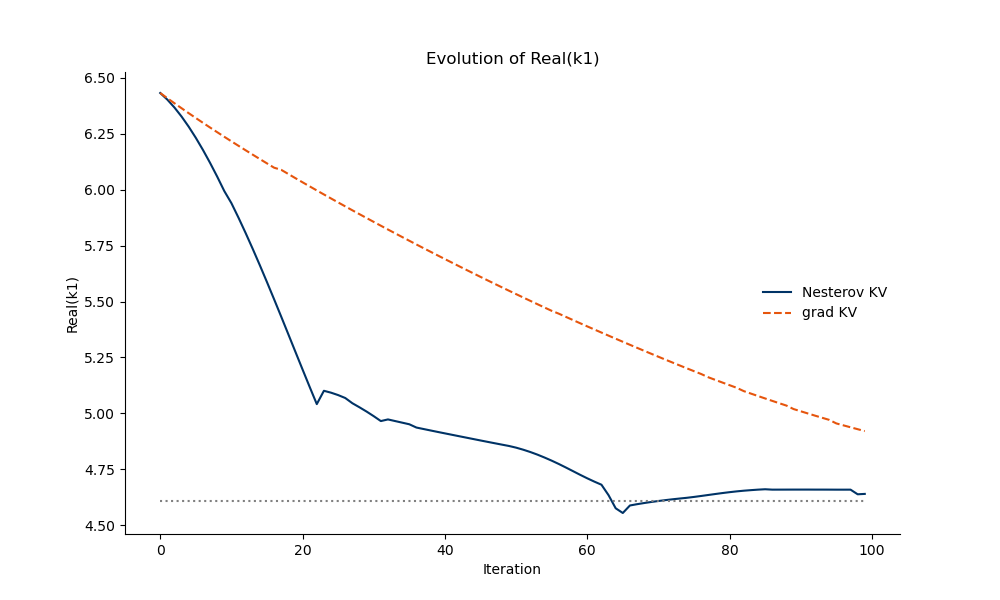}\label{fig:7b}
  }
  \caption{Comparison of Nesterov and gradient descent methods}
  \label{fig:comparisongrad}
\end{figure}

We here present some experiments to compare the behaviors of both methods in our case. Figure~\ref{fig:7a} presents the general comportment we can observe. The parameters used are same than the one used to obtain Figure~\ref{fig:robustexpected} but we consider the deterministic case here. Hence, generally, we do not observe a systematic performance gain for the Nesterov method in our specific framework. This can be explained by the nature of the cost functional. Moreover, the presence of restarting somehow stops the accelerating procedure we can find in Nesterov. Nevertheless, as shown in Figure \ref{fig:7b} where the reconstruction of the wavenumber only is done with the same parameters used to obtain Figure~\ref{fig:deterministparameter}, there are favorable cases where the Nesterov descent proves to be faster, justifying its use in our algorithm.

\subsubsection{The least squares functional}
In this article, we chose to study the Kohn-Vogelius functional instead of the usual least squares functional. In the case of the Laplacian, a comparison between those two functionals has been made in \cite{Afraites2007} showing that the Kohn-Vogelius functional is more robust than the least squares functional. We present here some numerical experiments that lead to the same conclusion in our case. In our context, the least squares functional writes:
$$\mathcal{J}(\Od_1,k_1) = \frac{1}{2} \int_{\Gamma_m} |u_{\rm R} - u_{\mathrm{obs}}|^2 ,$$ 
where $u_{\rm R} \in V$ solves~\eqref{eq:urud}. 
Using similar computations as for the Kohn-Vogelius functional, the derivatives of the least squares functional can be expressed as follows:

\[
    \DD_{\Od_1}\mathcal{J}(\Od_1,k_1)(\theta) = \mathfrak{Re}\left(\int_{\partial \Od_1} (k_1^2 - k_2^2) \rho u_{\rm R} (\theta \cdot \nn) \right) \quad \mbox{ and } \quad \frac{\partial \mathcal{J}(\Od_1,k_1)}{\partial k_1} = 2\int_{\partial \Od_1} \rho u_{\rm R} k_1,
\]

where $\rho \in V$ is the solution to the adjoint equation \begin{equation}\label{eq:adjointLS}\left\{ \begin{array}{rcll}
    \Delta \rho + k^2 \rho &=& 0 &\text { in } \mathcal{D}\\
    \partial_{\nn} \rho - ik\rho &=& 0 & \text{ on } \Gamma_{\rm c} \\
    \partial_{\nn} \rho - ik\rho &=& \overline{u_{\rm R} - u_{\mathrm{obs}}} & \text{ on } \Gamma_{\rm m} \\
    \left[\rho\right] \; = \; \left[ \partial_{\nn} \rho \right] & = & 0 & \text{ on }  \partial \Od_1 . \\
\end{array}
\right.
\end{equation}

\begin{figure}[h!tb]
  \centering
  \subfloat{
    \includegraphics[trim={1cm 0 2cm 1cm}, clip, width=0.48\textwidth]{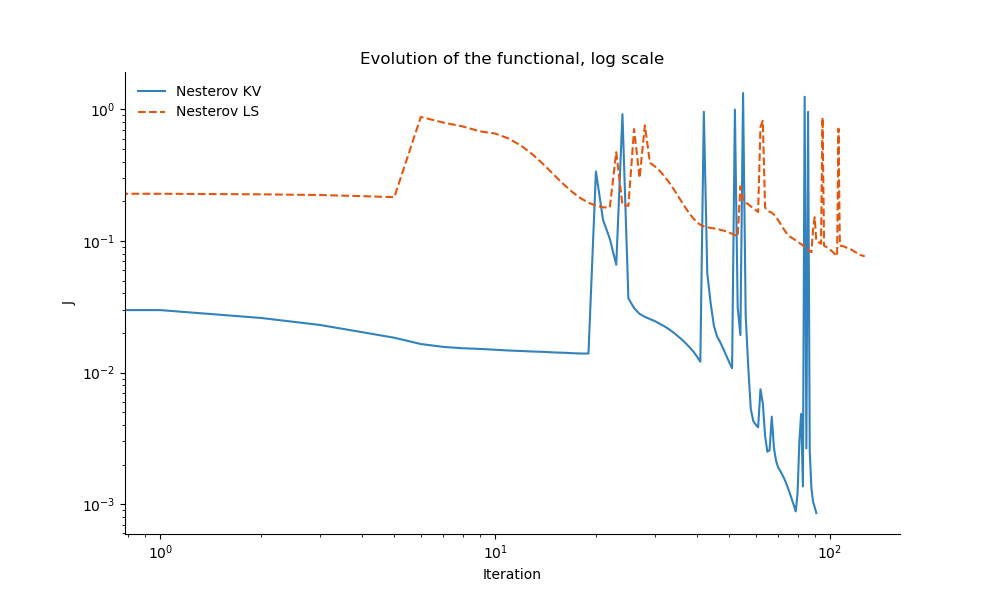}
  } 
  \hfill
  \subfloat{
    \includegraphics[trim={1.5cm 0 2.5cm 0cm}, clip, width=0.48\textwidth]{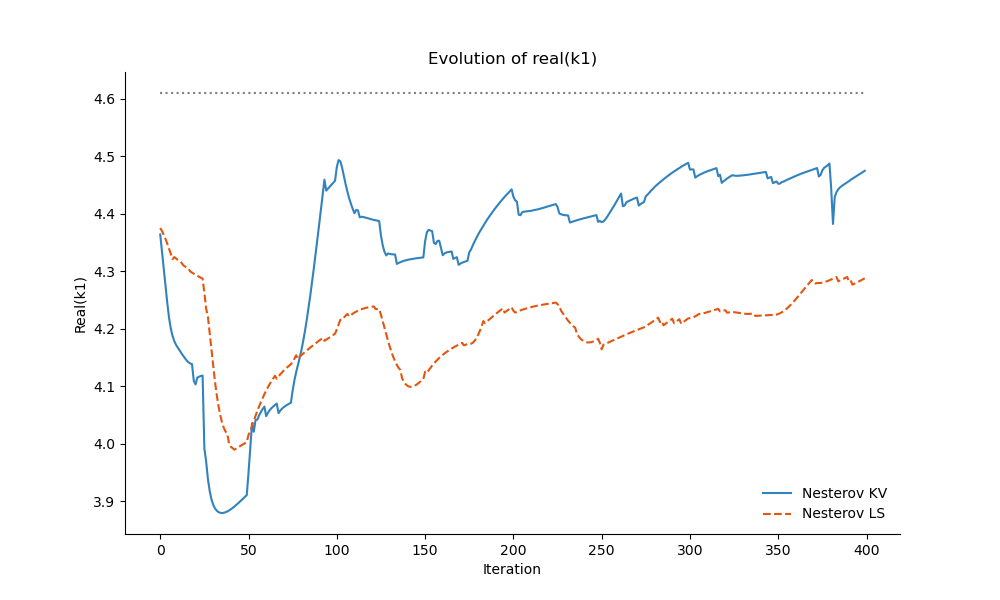}
  }
  \caption{Comparison of least squares and Kohn-Vogelius functional}
  \label{fig:comparison}
\end{figure}

We used the same parameter than in Figure \ref{fig:simultaneous}, but with a determinist data to perform the comparison.
We make several comments about Figure \ref{fig:comparison}. First, although the least squares functional is applicable to our inverse problem, its evolution during the reconstruction process reveals it to be "flatter" than the Kohn-Vogelius functional. This characteristic suggests a lack of robustness, a point previously noted by the authors in \cite{Afraites2007}. \\

In terms of computational cost, the Kohn-Vogelius functional requires solving twice as many equations to compute its derivatives. Nevertheless, this does not result in a doubling of the total execution time: the mesh is computed only once before those equations, that have all the same type. In practice, we observe that the Kohn-Vogelius functional requires fewer iterations to converge than the least squares functional, making this approach faster.

\section{Proofs}
\label{sec:proofs}
In this section, we provide the proofs of the propositions stated above. 

\subsection{Proof of the well-posedness of the problem \eqref{eq:transmission}}
\label{subsec:well-posed}

Since $u_{\rm obs} \in \HH^{1/2}(\Gamma_{\rm m})$, there exists a function $u_0 \in \HH^1(\D)$ such that the trace of $u_0$ coincides with $u_{\rm obs}$ on~$\Gamma_{\rm m}$. Up to studying the equation satisfied by $u - u_0$, we may assume in the following that $u_{\rm obs} = 0$ in the problem \eqref{eq:transmission}. Hence, all the problems considered in this paper are well-defined including the adjoint problems~\eqref{eq:adjoints}. We remind first two results of elliptic PDEs used in the following proof. 

 \begin{theorem}[Banach-Nečas-Babuska Theorem (BNB), {\cite[Theorem~25.9]{ErnGuermond2021}}]\label{theorem:BNB}
    Let $V$ be a Banach space, $W$ be a reflexive Banach space and $a \in \mathcal{L}(V \times W, \mathbb{C})$. The following items are equivalent:

    \begin{enumerate}
        \item For any $l \in W'$, there exists a unique $u \in V$ such that 
        \[a(u,w) = l(w), \text{ for all } w \in W.\]
        \item There exists $\beta > 0$ such that $\mathrm{inf}_{v \in V} \mathrm{sup}_{w \in W} \frac{a(v,w)}{||w||_W} = \beta$ and $\left(\forall v \in V, a(v,w) = 0 \Rightarrow (w = 0)\right).$
        \item There exists $\beta_1, \beta_2 > 0$ such that $\mathrm{inf}_{v \in V} \mathrm{sup}_{w \in W} \frac{a(v,w)}{||v||_V||w||_W} = \beta_1$ and $\mathrm{inf}_{w \in W} \mathrm{sup}_{v \in V} \frac{a(v,w)}{||v||_V||w||_W} = \beta_2.$
    \end{enumerate}

\end{theorem}

We also have a lemma that provides a sufficient condition for applying the BNB:

\begin{lemme}[Gårding's lemma, {\cite[Lemma~35.3]{ErnGuermond2021}}]\label{lem:garding}
Let $V \hookrightarrow L$ be two Banach spaces with compact embedding. Let $a: V \times V \rightarrow \mathbb{C}$ be a bounded sesquilinear form. Assume that there exist two real numbers $\beta, \gamma > 0$ such that the following holds:

\begin{center}
\begin{enumerate}
    \item $|a(v,v)| + \beta||v||_L^2 \ge \gamma ||v||_V^2, \forall v \in V$
    \item $\left[a(v,w) = 0, \forall w \in V \right] \Rightarrow (v = 0).$
\end{enumerate}
\end{center}

Then the second condition of the BNB Theorem is satisfied.
\end{lemme}

\begin{proof}[Proof of the proposition~\ref{prop:wpforward}]
To prove this we make use of the Lemma~\ref{lem:garding} and thus prove its two points. Using the notations of this lemma, let $V = \{u \in \HH^1(\D) \text{ such that } [u] = 0 \text{ on } \Od_1 \text{ and } u = 0 \text{ on } \Gamma_{\rm m}\}$ and $L = \mathrm{L}^2(\Od_1)$.
\begin{enumerate}
        \item 
    Let us denote $A(u,v) = \int_{\D} \nabla u \nabla \overline{v} - k^2 u \overline{v} - ik_2 \int_{\Gamma_{\rm c}} u \overline{v}$ and consider $\mathfrak{Re}(A(u,-ik_2u))$.
    \begin{align*}
        & \mathfrak{Re}(A(u,-ik_2u)) \\
        &= \mathfrak{Re}\left(i \left[\overline{k_2} ||\nabla u||_{\mathrm{L}^2(\Od_1)}^2 + \overline{k_2} ||\nabla u||_{\mathrm{L}^2(\Od_2)}^2 - k_1^2 \overline{k_2} ||u||_{\mathrm{L}^2(\Od_1)}^2 - k_2^2 \overline{k_2} ||u||_{\mathrm{L}^2(\Od_2)}^2 \right]+ |k_2|^2 ||u||_{\mathrm{L}^2(\Gamma_{\rm c})}^2 \right) \\
        &= \mathfrak{Im}(k_1^2 \overline{k_2})||u||_{\mathrm{L}^2(\Od_1)}^2 + \mathfrak{Im}(k_2) \left[ ||\nabla u||_{\mathrm{L}^2(\D)}^2 + ||k_2u||_{\mathrm{L}^2(\Od_2)}^2 \right] +  |k_2|^2 ||u||_{\mathrm{L}^2(\partial \D)}^2. 
    \end{align*}
    Hence, with the definition of $\mathfrak{K}_{\rm adm}$, $$|k_2||A(u,u)| + |\mathfrak{Im}(k_1^2 \overline{k_2})| ||u||_{\mathrm{L}^2(\Od_1)}^2\ge \mathfrak{Im}(k_2) \left[ ||\nabla u||_{\mathrm{L}^2(\D)}^2 + ||k_2u||_{\mathrm{L}^2(\Od_2)}^2 \right]  + |k_2|^2 ||u||_{\mathrm{L}^2(\Gamma_{\rm c})}^2$$ and we get the first point of the Gårding's lemma. 

    \item For the second point, assume that $A(u,v) = 0$ for every $v \in V$. Hence, 
    \[\mathfrak{Im}(A(u,u)) = - \mathfrak{Re}(k_2) \left( ||u||_{\mathrm{L}^2(\Gamma_{\rm c})} + 2 \mathfrak{Im}(k_2) ||u||_{\mathrm{L}^2(\Od_2)} \right) - 2 \mathfrak{Re}(k_1) \mathfrak{Im}(k_1) ||u||_{\mathrm{L}^2(\Od_1)} = 0\] and, since $\mathfrak{Im}(k_1) > 0$ and the real parts $\mathfrak{Re}(k_1)$ and $\mathfrak{Re}(k_2)$ have the same sign, we obtain $||u||_{\mathrm{L}^2( \Od_1 )} = ||u||_{\mathrm{L}^2( \Od_2 )} =  0$. 
    \end{enumerate}
    
\end{proof}

\subsection{Proof of the uniqueness}

The approach is usual concerning the identifiability for inverse problems. First, we write an orthogonality formula \eqref{eq:orthogonal} and then we use the construction of periodical special solutions to the Helmholtz equations built in \cite{hahner}.

\begin{proof}[Proof of the proposition \ref{prop:uniqueness}]
Let $u_{\rm a} \in V_{\rm a}$ and $u_{\rm b} \in V_{\rm b}$ be the solution $u_{\RR}$ of~\eqref{eq:urud} with $\Od_1$, $k_1$ and $h$ being respectively~$\Od_{\rm a}, k_{\rm a}, h_{\rm a}$ and $\Od_{\rm b}, k_{\rm b}, h_{\rm b}$, where $V_{\rm a}$ and $V_{\rm b}$ are the corresponding variational spaces. We note~$\Lambda$ their common Robin-to-Dirichlet operator given in~\eqref{def:RobToDirOp}. Finally, we define the function $\delta$, as being the difference of the two functions, constant by parts, $k_{\rm a}^2$ and $k_{\rm b}^2$. Notice that $\delta = 0$ in $\D \setminus \Od_{\rm a} \cup \Od_{\rm b}.$ Notice also, that $u_{\rm b} \in V_{\rm a}$ and $u_{\rm a} \in V_{\rm b}$, by inner elliptic regularity.\\

Writing the variational formulation defining $u_{\rm a}$ and applying it with $v = \overline{u_{\rm b}} \in V_{\rm a}$, we get 
$$\int_{\D} \nabla u_{\rm a} \nabla u_{\rm b} - k^2 u_{\rm a} u_{\rm b} - ik_2 \int_{\partial D} u_{\rm a} u_{\rm b} = \int_{\partial D} h_{\rm a} u_{\rm b}.$$
Similarly, we compute the variational form defining $u_{\rm b}$ with $\overline{u_{\rm a}} \in V_{\rm b}$ as a test function. Computing the difference of both expressions, we have $$\int_{\D} \delta u_{\rm a} u_{\rm b} = \int_{\partial \D} h_{\rm b} u_{\rm a} - h_{\rm a} u_{\rm b} = \int_{\partial \D} h_{\rm b} \Lambda(h_{\rm a}) - h_{\rm a} \Lambda(h_{\rm b}).$$ With the Green formula, we show the right term of this equation is $0$ and then 
\begin{equation}\label{eq:orthogonal}
    \int_{\D} \delta u_{\rm a} u_{\rm b} = 0.
\end{equation}

Let $B$ a ball such that $\D \subset B$ and consider $q: B \rightarrow \mathbb{C}$ a continuous function with compact support. Let~$\xi \in \mathbb{C}^d$ non-zero such that $\xi \cdot \xi = 0$. We are looking for specials solutions to the Helmholtz equation 
\begin{equation}\label{eq:special}
\Delta u + q u = 0 
\end{equation}
of the form $u_{\xi}(x) = \mathrm{e}^{i\xi \cdot x}(1+r_{\xi}(x))$. We remark that if we take $$q_{\rm a} = \begin{cases} k_{\rm a}^2 \in \Od_{\rm a} \\ k_2^2 \in \D \setminus \Od_{\rm a} \\ 0 \in B\setminus \D,
\end{cases}$$ then a solution to the Helmholtz equation \eqref{eq:special} is a solution to the problem $\Tau(\Od_{\rm a}, k_{\rm a})$ in the domain $\D$. Of course, this can be done if we replace ${\rm a}$ by ${\rm b}$. Using \cite[Lemma~5]{hahner}, this solution $u_{\xi}$ exists and satisfies for some constant $C > 0$, independent of $\xi$, 
$$||r||_{\mathrm{L}^{\infty}(B)} \le \frac{C}{|\boldsymbol{\mathfrak{Im}}(\xi)|}.
$$
From this, since the dimension $d$ is $3$, we follow the construction done in~\cite[Theorem 6]{hahner}: there exist~$u_{\rm a}(t)$ and~$u_{\rm b}(t)$ special solutions such that $r_{\rm a}(t)$ and $r_{\rm b}(t)$ go to $0$ when $t$ goes to $\infty$. Hence, with \eqref{eq:orthogonal},
$$
\int_{\D} \delta \mathrm{e}^{i y \cdot x}(1+r_{\rm a}(t))(1 + r_{\rm b}(t)) = 0, \text{ for all } y \in \mathbb{C}^{d}, \quad \text{ for all } t \ge 0.
$$ 
Finally, letting $t$ goes to $+ \infty$, we deduce that $\delta = 0$ by Fourier inversion. This means $k_{\rm a} = k_{\rm b}$ and also~$\Od_{\rm a} = \Od_{\rm b}$.
\end{proof}
\subsection{Reminder on shape optimization}
First, we remind two important lemmas that are necessary for the following proof.
\begin{lemme} \label{lem:jprime} Let $J: \mathcal{D}_{adm} \rightarrow \mathbb{R}$ be a functional defined on a class of admissible shapes $\D_{adm}$. Let $\theta \in \mathrm{W}_0^{1,\infty}(\D)$ and~$f \in \mathrm{L}^1(\mathbb{R}^d)$. 
If $J(\Od_{\theta}) = \int_{\Od_\theta} f(\theta)$, then $$J'(\Od)(\theta) = \int_{\Od} f'(x) + \int_{\partial \Od} f (\theta \cdot \nn),$$ where $f'(x)$ denotes the Fréchet derivative of the map $\theta \mapsto f(x + \theta(x))$.
\end{lemme}
\begin{proof}
    This is the result \cite[corollary 5.2.8]{henrotpierre2018} from the book of Henrot and Pierre applied with $\Phi(x) = (I + t \theta(x)).$
\end{proof}

\begin{lemme}[{\cite[Lemma 5.3.3]{henrotpierre2018}}]\label{lem:shapederivative}
    Let $\psi: \mathrm{W}^{1,\infty}(\bR^d) \rightarrow \mathrm{W}^{1,\infty}(\bR^d)$ continuous at $0$ with $\psi(0) = I$ and $g: \mathrm{W}^{1,\infty}(\bR^d) \rightarrow \mathrm{W}^{1,p}(\mathbb{R}^d)$, with $1 \le p < + \infty$. If $\theta \mapsto (g(\theta),\psi(\theta)) \in \mathrm{L}^p(\bR^d) \times \mathrm{L}^{\infty}(\bR^d)$ is differentiable at $0$ with $g'(0)$ continuous, then the application 

    \[G: \theta \in \mathrm{W}^{1,\infty}(\bR^d) \mapsto g(\theta) \circ \psi(\theta) \in \mathrm{L}^p(\mathbb{R}^d)\]

 is differentiable at $0$ and

    \[G'(0)\xi = g'(0)\xi + \nabla g(0) \cdot \psi'(0) \xi.\]
\end{lemme}

\subsection{Proofs concerning the deterministic case}
The proof of the proposition~\ref{prop:adjointfirststrategy} is divided into two steps. First, we classically differentiate the functional $\mathcal{K}$ using the technical lemmas recalled in the last subsection. We obtain an expression which is not satisfactory from a numerical point of view. Indeed, in order to apply a Nesterov descent scheme, we need to efficiently compute a clear direction of descent. Thus, the second step of the proof consists in introducing adjoints in order to establish the proposition \ref{prop:adjointfirststrategy}.

\begin{lemme}\label{prop:derivativefirststrategy}
    For every $\Od$ admissible, and $\theta \in \WW_0^{1,\infty}(\D)$ the shape derivative of the functional \eqref{def:K} exists and

    \begin{equation}\label{Kprime}
    \DD_{\Od_1}\mathcal{K}(\Od_1,k_1)(\theta) = \mathfrak{Re} \left(\int_{\D} \nabla (u_{\DD}'-u_{\RR}') \cdot \nabla (\overline{u_{\DD}-u_{\RR}}) \right),
    \end{equation}
where $u_{\RR}' \in \HH^1(\D)$ and $u_{\DD}' \in \HH^1(\D)$ are solutions to
\[
\begin{array}{cc}
\left\{
    \begin{array}{rcll}
    \Delta u_{\RR}' + k^2 u_{\RR}' & = & 0 &  \text{ in }  \D\\
    \partial_{\nn} u_{\RR}' - iku_{\RR}' & = & 0 & \text{ on }  \partial \D \\
    \left[u_{\RR}'\right] & = & 0 & \text{ on }  \partial \Od_1 \\
    \left[\partial_{\nn} u_{\RR}'\right] & = & (k_1^2 - k_2^2) u_{\RR} (\theta \cdot \nn) & \text{ on }  \partial \Od_1, \\

\end{array}\right.
&
\left\{
\begin{array}{rcll}
    \Delta u_{\DD}' + k^2 u_{\DD}' & = & 0 & \text{ in } \D \\
    \partial_{\nn} u_{\DD}' - iku_{\DD}' & = & 0 &  \text{ on } \Gamma_{\mathrm{c}} \\
    u_{\DD}' & = & 0 & \text{ on } \Gamma_{\mathrm{m}} \\
    \left[u_{\DD}'\right] & = & 0 & \text{ on }  \partial \Od_1 \\
    \left[ \partial_{\nn} u_{\DD}'\right] & = & (k_1^2-k_2^2)u_{\DD} (\theta \cdot \nn) & \text{ on } \partial \Od_1.
\end{array}
\right.
\end{array}\]
\end{lemme}

 The proof of this lemma is similar to the one proposed by Afraites and Dambrine in~\cite[Theorem 3.1]{Afraites2007} for a different problem. However, our weak formulation differs and we have to proceed differently for the steps~$2$ and~$3$ of their proof. \\

\begin{proof}[Proof of the lemma~\ref{prop:derivativefirststrategy}]
Fix $\Od \in \mathcal{U}_{\mathrm{adm}}$ and $\theta \in \mathrm{W}_0^{1,\infty}(\mathcal{D})$. We denote by $u_{\DD}(\theta)$ and $u_{\RR}(\theta)$ the solutions to the problems \eqref{eq:urud} with $\Od_1 = (I + \theta)(\Od)$. We only need to prove that $\theta \mapsto u_{\RR}(\theta)$ and $\theta \mapsto u_{\DD}(\theta)$ have shape derivatives. The two computations being very similar, we will do this work for $u_{\RR}$ only. To lighten the notations, we denote $u_\theta= u_{\RR}(\theta)$. We can notice that for every $\theta \in \mathrm{W}_0^{1,\infty}(\D)$, the functions $u_\theta$ are in~$\mathrm{H}_{\Od_\theta}^1(\D)$. To study the derivation of the application $\theta \mapsto u_\theta$, we consider the application $w(t) = u_\theta \circ (I + t\theta)$ defined on the space $\mathrm{H}^1_{\Od}(\D)$, which does not depend on the choice of $\theta$. \\

 If $u$ has a derivative $\dot{u}$, then we should have $\frac{w(t)-u}{t}$ that goes to $\dot{u}$. This remark gives us a candidate for the variational formulation that should be satisfied by $\dot{u}$. We will use that formulation to show that $\frac{w(t)-u}{t}$ effectively goes to $\dot{u}$ and then compute the equation satisfied by the shape derivative $u'$ thanks to the lemma \ref{lem:shapederivative}. \newline

 By a change of variable $w(t)$ is solution of the weak problem 

\[\int_{\D} \left(A(t) \nabla w \nabla \overline{v} - J(t) k^2 w \overline{v}\right) - ik_2 \int_{\partial \D} w \overline{v} = \int_{\partial \D} h \overline{v}, \quad \forall v \in V\]

 where $J(t) = \mathrm{det}(I + tD\theta)$ and $A(t) = J(t)(I + tD\theta)^{-1}(I + tD\theta^T)^{-1}$. Subtracting by the variational form associated to $u$ and dividing by $t \neq 0$, we get

\[\int_{\D} \left( A(t)\frac{\nabla w(t) - \nabla u}{t} \nabla\overline{v} - J(t) k^2 \frac{w(t) - u}{t} \overline{v}\right) - ik_2 \int_{\partial \D} \frac{w(t)-u}{t} \overline{v} = \int_{\D} \frac{I-A}{t} \nabla u \nabla \overline{v} - k^2 \frac{1-J}{t} u \overline{v}, \forall v \in V\]

 The matrix $A$ is usual in shape optimization and satisfies the following properties \cite{Afraites2007}: $A(t,x)$ is symmetric positive, $A(0) = I$ with $$\mathcal{A} = \frac{\mathrm{d}}{\mathrm{dt}}\restriction{A(t)}{t=0} = \mathrm{div}(\theta) I - (D\theta^T + D\theta)$$ and for $t$ small enough, we have \begin{equation}\label{eq:A}
    \forall X \in \bR^d, X^T A(t) X \ge \frac{||X||^2}{2}.
\end{equation}Similarly, the quotient $\frac{J(t)-1}{t}$ goes to $\mathrm{div}(\theta)$ when $t$ goes to $0$ by differentiation of the $\mathrm{det}$ application. Therefore, we let $\dot{u} \in V$ that satisfies

\[\int_{\D} \left(\nabla \dot{u} \nabla \overline{v} - k^2 \dot{u} \overline{v} \right)- ik_2 \int_{\partial \D} \dot{u} \overline{v} = -\int_{\D} \mathcal{A} \nabla u \nabla \overline{v} - k^2 \mathrm{div}(\theta) u \overline{v}, \quad \forall v \in V.\]

 The function $\dot{u}$ is well defined since this is the same variational form than the transmission problem, but with a different sesqui-linear term on the right. At this stage, we do not know if $\dot{u}$ is the desired derivative. It is merely a candidate function that exists. Denote $\frac{w(t)-u}{t} - \dot{u}$ by $d(t)$. We are done if we show that $d(t)$ goes to $0$ in $\mathrm{H}^1(\D)$. Subtracting the two last equations, we get 

\begin{align*}
    \int_{\D} A(t) \nabla d(t) \nabla \overline{v} - k^2 J(t) d(t) \overline{v} - ik_2 \int_{\partial \D} d(t) \overline{v} =& \int_{\D} \left[ \left(\frac{I-A(t)}{t}+\mathcal{A}\right)\nabla u + (I-A(t))\dot{u}\right] \nabla \overline{v}\\
    &- \int_{\D} k^2 \left[\left(\frac{1-J(t)}{t}+ \mathrm{div}(\theta)\right)u + (1-J(t)) \dot{u}\right]\overline{v}.
\end{align*}

 We use $d(t)$ as a test function and proceed as in the last proof, by majoring the imaginary part of the left side by the modulus of the right side. Thus, using similar method than in the sub-section \ref{subsec:well-posed}, we get by the Cauchy-Schwarz inequality and the continuity of $J$, for $t$ small enough, 

\[||d(t)||_{\mathrm{L}^2(\partial \D)}^2 \lesssim_{u,k,\dot{u}} \left(\left| \frac{I-A(t)}{t}+\mathcal{A} \right|_{\infty} + \left|I-A(t)\right|_{\infty} + \left|\frac{1-J(t)}{t}+\mathrm{div}(\theta)\right| + |1-J(t)|\right)||d(t)||_{\mathrm{H}^1(\D)}\]
 and 
\[||d(t)||_{\mathrm{L}^2(\D)}^2  \lesssim_{u,k,\dot{u}} \left(\left| \frac{I-A(t)}{t}+\mathcal{A} \right|_{\infty} + \left|I-A(t)\right|_{\infty} + \left|\frac{1-J(t)}{t}+\mathrm{div}(\theta)\right| + |1-J(t)|\right) ||d(t)||_{\mathrm{H}^1(\D)},\] where $\lesssim_{A,B}$ denotes being less than or equal to, up to a multiplicative constant that depends only on the quantities $A$ and $B$. \\

Finally, if we come back to the whole expression with $d(t)$ as a test function, by the property \eqref{eq:A} of $A$, we get
\[||d(t)||_{\mathrm{H}^1(\D)}^2  \lesssim_{u,k,\dot{u}} \left(\left| \frac{I-A(t)}{t}+\mathcal{A} \right|_{\infty} + \left|I-A(t)\right|_{\infty} + \left|\frac{1-J(t)}{t}+\mathrm{div}(\theta)\right| + |1-J(t)|\right)||d(t)||_{\mathrm{H}^1(\D)} \]
 so that $||d(t)||_{\mathrm{H}^1(\D)} \rightarrow 0$ when $t$ goes to zero. We have just shown that $\dot{u}$ is the derivative of $u$. Consider now the shape derivative $u' = \dot{u} - \nabla u \cdot \theta$ (given by the lemma \ref{lem:shapederivative} with $\psi(\theta) = (I + \theta)^{-1}$ and $g(\theta) = u_\theta \circ (I + \theta)$). Adapting the computation of the step 4 of the proof of \cite[Theorem~3.8]{Afraites2007}, we first conclude that

\[\int_{\D} \nabla (\dot{u} - \nabla u \cdot \theta) \nabla \overline{v} = \int_{\D} \left(\nabla u' \nabla \overline{v} + \mathcal{A} \nabla u \nabla \overline{v} - k^2 u (\theta \cdot \nabla \overline{v})\right).\]

 Using the variational formulation satisfied by $\dot{u}$, we get
\begin{align*}
    \int_{\D} \left(\nabla u' \nabla \overline{v} - k^2 u' \overline{v} \right)- ik_2 \int_{\partial \D} u' \overline{v} &= \int_{\D} k^2 u (\theta \cdot \overline{v}) + \left[(\nabla u \cdot \theta)+ k^2 \left( \mathrm{div}(\theta) u  + (\nabla u \cdot \theta) \right)\right] \overline{v}\\
    &= \int_{\D} k^2 u (\theta \cdot \nabla \overline{v}) + k^2 \mathrm{div}(u \theta) \overline{v} \\
    &= \int_{\D} k^2 \mathrm{div}(u \theta \overline{v}).
\end{align*}

 By the Green formula on both $\Od_1$ and $\Od_2$, we have

\[\int_{\D} \left(\nabla u' \nabla \overline{v} - k^2 u' \overline{v} \right)- ik_2 \int_{\partial \D} u' \overline{v} = \int_{\partial \Od_2 \cap \Od_1} k_2^2 u \overline{v} (\theta \cdot \nn) + \int_{\partial \Od_1} k_1^2 u\overline{v} (\theta \cdot \nn) = \int_{\partial \Od_1} (k_1^2 - k_2^2) u \overline{v} (\theta \cdot \nn).\]

 We recognize the weak formulation of the problem \eqref{eq:transmission}. By the lemma $\ref{lem:jprime}$ below, we conclude that shape derivative exists and has the desired expression.
\end{proof} 

\begin{proof}[Proof of the adjoint method: proposition \ref{prop:adjointfirststrategy}]

By an integrations by part and the equation satisfied by $\rho$ we have:
    \begin{align*}
        \int_{\Od_i} \nabla u_{\DD}' \cdot \nabla (\overline{u_{\DD}-u_{\RR}}) &= \int_{\Od_i} -(u_{\DD}' \Delta (\overline{u_{\DD}-u_{\RR}}) + \int_{\partial \Od_i} u_{\DD}' \partial_{\nn} (\overline{u_{\DD}-u_{\RR}}) \\
        &= \int_{\Od_i} - u_{\DD}' (\Delta \rho_{\DD} + k^2 \rho_{\DD}) + \int_{\partial \Od_i} u_{\DD}' \partial_{\nn} (\overline{u_{\DD} - u_{\RR}}) \\
        &= \int_{\partial \Od_i} \partial_{\nn} u_{\DD}'\rho_{\DD} - \partial_{\nn} \rho_{\DD} u_{\DD}' + u_{\DD}' \partial_{\nn} (\overline{u_{\DD} - u_{\RR}}).
    \end{align*}

Similarly,
    \[\int_{\Od_i} \nabla u_{\RR}' \cdot \nabla (\overline{u_{\DD}-u_{\RR}}) = \int_{\partial \Od_i} \partial_{\nn} u_{\RR}' \rho_{\RR} + u_{\RR}' ( \partial_{\nn} (\overline{u_{\DD}-u_{\RR}}) - \partial_{\nn} \rho_{\RR}).\]
On $\partial \D$, the integrals are null thanks to the adjoint equation. Hence, we just need to compute the quantities
\[\mathcal{I}_{\DD} =\int_{\partial \Od_1} \partial_{\nn} u_{\DD}' \rho_{\DD} + u_{\DD}' ( \partial_{\nn} (\overline{u_{\DD}-u_{\RR}}) - \partial_{\nn} \rho_{\DD}) + \int_{\partial \Od_2 \cap \Od_1} \partial_{\nn} u_{\DD}' \rho_{\DD} + u_{\DD}' ( \partial_{\nn} (\overline{u_{\DD}-u_{\RR}}) - \partial_{\nn} \rho_{\DD})\] and
\[\mathcal{I}_{\RR} =\int_{\partial \Od_1} \partial_{\nn} u_{\RR}' \rho_{\RR} + u_{\RR}' ( \partial_{\nn} (\overline{u_{\DD}-u_{\RR}}) - \partial_{\nn} \rho_{\RR})+ \int_{\partial \Od_2 \cap \Od_1} \partial_{\nn} u_{\RR}' \rho_{\RR} + u_{\RR}' ( \partial_{\nn} (\overline{u_{\DD}-u_{\RR}}) - \partial_{\nn} \rho_{\RR}).\]

 By the jump conditions and the identity $[ab] = [a] \frac{b_1 + b_2}{2} + [b] \frac{a_1 + a_2}{2}$, we have
\[\mathcal{I}_{\DD} = \int_{\partial \Od_1} (k_1^2-k_2^2) u_{\DD} \rho_{\DD}(\theta \cdot \nn) \quad \text{and} \quad \mathcal{I}_{\RR} = \int_{\partial \Od_1} (k_1^2-k_2^2) u_{\RR} \rho_{\RR} (\theta \cdot \nn)\] so that
\[\mathrm{D}_{\Od_1}\mathcal{K}(\Od_1,k_1)(\theta) = \mathfrak{Re} \left( \int_{\D} \nabla (u_{\DD}' - u_{\RR}') \cdot \nabla \overline{u_{\DD} - u_{\RR}}\right) = \mathfrak{Re}\left(\int_{\Od_1} (k_1^2 - k_2^2) (u_{\DD} \rho_{\DD} - u_{\RR} \rho_{\RR}) (\theta \cdot \nn)\right).\]
\end{proof}

To study the derivative with respect to the wavenumber, we proceed in a similar way. First, we compute a preliminary expression of the derivative. Then, we introduce the adjoint variables $\rho_D$ and $\rho_R$ in order to prove the proposition~\ref{prop:derivativekfirststrategy}. 
\begin{lemme}[Derivative with respect to $k$]\label{lemma:derivativekfirststrategy}
For every $k_1$ admissible, the derivative of the functional \eqref{def:K} with respect to $k_1$ exists and

    \begin{equation}\label{Kprimek1}
    \frac{\partial \mathcal{K}(\Od_1,k_1)}{\partial k_1} =\int_{\D} \nabla (u_{\DD}-u_{\RR}) \cdot \nabla (\overline{u_{\DD}'-u_{\RR}'}),
    \end{equation}
where $u_{\RR}'$ and $u_{\DD}' \in \HH^1(\D)$ are solution to
\[
\begin{array}{cc}
\left\{
    \begin{array}{rcll}
    \Delta u_{\RR}' + k_2^2 u_{\RR}' & = & 0 & \text{ in } \Od_2\\
    \Delta u_{\RR}' + k_1^2 u_{\RR}' & = & -2k_1u_{\RR} & \text{ in }  \Od_1\\
    \partial_{\nn} u_{\RR}' - iku_{\RR}' & = & 0 & \text{ on } \partial \D \\
    \left[\partial_{\nn} u_{\RR}'\right] \; = \; \left[u_{\RR}'\right] & = & 0 & \text{ on }  \partial \Od_1, \\
\end{array}\right.
&
\left\{
\begin{array}{rcll}
    \Delta u_{\DD}' + k^2 u_{\DD}' & = & 0 & \text{ in } \Od_2 \\
    \Delta u_{\DD}' + k_1^2 u_{\DD}' & = & -2k_1u_{\DD} & \text{ in } \Od_1\\
    \partial_{\nn} u_{\DD}' - iku_{\DD}' & = & 0 & \text{ on } \Gamma_{\mathrm{c}} \\
    u_{\DD}' & = & 0 & \text{ on } \Gamma_{\mathrm{m}} \\
    \left[ \partial_{\nn} u_{\DD}'\right] \; = \; \left[u_{\DD}'\right] & = & 0 & \text{ on }  \partial \Od_1. \\
\end{array}
\right.
\end{array}\]
\end{lemme}

 Contrary to the proof of proposition \ref{prop:derivativefirststrategy}, we will use the implicit function theorem, that is more convenient in our case. As before, we will only consider the derivatives of $k_1 \mapsto u_{\RR}(k_1)$. \begin{proof} Let us fix an admissible wave-number $k_1$ and define the function \[F: \begin{cases}
        \bC \times\HH^1(\D) \rightarrow (H^1(\D))' \\
        (k,u) \mapsto \left(v \mapsto \int_{\D} \nabla u \nabla \overline{v} - k^2 u \overline{v} - ik\int_{\partial \D} (u + h) \overline{v}\right).
    \end{cases}\] The application $F$ is $\mathcal{C}^1$ and is chosen such that $F(k,u) = 0$ if and only if $u$ is solution to the problem~\eqref{eq:transmission}. We have $$\partial_u F(k,u_{\RR}(k))(u)v = \int_{\D} \nabla u \nabla \overline{v} -k^2 u \overline{v} - ik \int_{\partial \D} u \overline{v}.$$ This is an isomorphism from $\mathrm{H}^1(\D)$ to $(H^1(\D))'$  by well-posedness of the problem \eqref{eq:transmission}. By the implicit function theorem, there exists $\Tilde{u}: \C^0(\D,\bC) \rightarrow\HH^1(\D)$ a $\C^1$ function defined on a neighborhood of $0$ such that $F(k,\Tilde{u}(k)) = 0$. By well-posedness of the problem \eqref{eq:transmission}, we have $\Tilde{u}(\Od_1,k) = u_{\RR}(k)$ and $k \mapsto u_{\RR}(k)$ is $\C^1$ on a neighborhood of $k_1$. Then, derivating the variational form, we see that $u_{\RR}'(k_1)$ satisfy $$\int_{\D} \nabla u_{\RR}' \nabla \overline{v} - k^2 u_{\RR}' \overline{v} - ik_2 \int_{\partial \D} u_{\RR}' \overline{v} = \int_{\Od_1} 2k_1 u_{\RR} \overline{v}, \quad \forall v \in V.$$ This shows that $u_{\RR}'$ satisfy the desired problem. 
\end{proof}
The expression of the derivative presented in Lemma~\ref{prop:derivativekfirststrategy} is sufficient to compute a gradient efficiently. However, this approach requires us to compute $u_{\RR}'$ and $u_{\DD}'$ at each step. A more efficient alternative is to compute the gradient using the adjoints $\rho_{\DD}$ ans $\rho_{\RR}$ which have already been calculated. Consequently, the only additional cost is integrating a function over the domain to obtain the gradient with respect to k.
\begin{proof}[Proof of the proposition \ref{prop:derivativekfirststrategy}]

    We decompose the derivative as $$
    \frac{\partial \mathcal{K}(\Od_1,k_1)}{\partial k_1} = \mathfrak{Re}\left(\int_{\D} \nabla u_{\DD}' \cdot \nabla (\overline{u_{\DD}-u_{\RR}})\right) - \mathfrak{Re}\left( \int_{\D} \nabla u_{\RR}' \cdot \nabla (\overline{u_{\DD}-u_{\RR}})\right).$$ 
    By integration by parts and the equation satisfied by $\rho$ we have:
    \begin{align*}
        \int_{\D} \nabla u_{\DD}' \cdot \nabla (\overline{u_{\DD}-u_{\RR}}) &= \int_{\D} -(u_{\DD}' \Delta (\overline{u_{\DD}-u_{\RR}}) + \int_{\partial \D} u_{\DD}' \partial_{\nn} (\overline{u_{\DD}-u_{\RR}}) \\
        &= \int_{\D} - u_{\DD}' (\Delta \rho_{\DD} + k^2 \rho_{\DD}) + \int_{\partial \D} u_{\DD}' \partial_{\nn} (\overline{u_{\DD} - u_{\RR}}) \\
        &= -\int_{\D} \rho_{\DD} \left(\Delta u_{\DD}' + k^2 u_{\DD}'\right) + \int_{\partial \D} \partial_{\nn} u_{\DD}'\rho_{\DD} - \partial_{\nn} \rho_D u_{\DD}' + u_{\DD}' \partial_{\nn} (\overline{u_{\DD} - u_{\RR}})\\
        &= 2 \int_{\Od_1} \rho_{\DD} u_{\DD} k_1 + \int_{\partial \D} \partial_{\nn} u_{\DD}'\rho_D - \partial_{\nn} \rho_D u_{\DD}' + u_{\DD}' \partial_{\nn} (\overline{u_{\DD} - u_{\RR}}).
    \end{align*}

 Decomposing the boundary as $\partial \D = \overline{\Gamma_{\rm c} \cup \Gamma_{\mathrm{m}}}$, we study the integral on each of the three parts. On $\Gamma_{\mathrm{m}}$, $u_{\DD}' = \rho_{\DD} =0$ and the integral is zero. On $\Gamma_{\mathrm{c}}$, we write 

    \[\partial_{\nn} u_{\DD}' \rho_{\DD} - \partial_{\nn} \rho_{\DD} u_{\DD}' + u_{\DD}' \partial_{\nn}(\overline{u_{\DD}-u_{\RR}}) = u_{\DD}' (\partial_{\nn} \rho_{\DD} - ik \rho_{\DD} - \partial_{\nn} \rho_{\DD}) + \partial_{\nn} u_{\DD}' \rho_{\DD} = \rho_{\DD} ( \partial_{\nn} u_{\DD}' - ik u_{\DD}') = 0\]

 so that the integral is also 0. Finally, we get $\int_{\D} \nabla u_{\DD}' \cdot \nabla (\overline{u_{\DD}-u_{\RR}}) = 2 \int_{\Od_1} \rho_{\DD} u_{\DD} k_1$. Let us consider now the second term. The first step is exactly the same and we get

    \[
    \int_{\D} \nabla u_{\RR}' \cdot \nabla (\overline{u_{\DD}-u_{\RR}})= 2 \int_{\Od_1} \rho_{\RR} u_{\RR} k_1 + \int_{\partial \D} \partial_{\nn} u_{\RR}'\rho_{\RR} - \partial_{\nn} \rho_{\RR} u_{\RR}' + u_{\RR}' \partial_{\nn} (\overline{u_{\DD} - u_{\RR}}).
    \]

 On $\partial \D$, the terms that we integrate become

    \[\partial_{\nn} u_{\RR}' \rho_{\RR} - \partial_{\nn} \rho_{\RR} u_{\RR}' + u_{\RR}' \partial_{\nn}(\overline{u_{\DD}-u_{\RR}}) = u_{\RR}' (\partial_{\nn} \rho_{\RR} - ik \rho_{\RR} - \partial_{\nn} \rho_{\RR}) + \partial_{\nn} u_{\RR}' \rho_{\RR} = \rho_{\RR} ( \partial_{\nn} u_{\RR}' - ik u_{\RR}') = 0.\]

 Finally, we get $\int_{\D} \nabla u_{\RR}' \cdot \nabla (\overline{u_{\DD}-u_{\RR}}) = 2 \int_{\Od_1} \rho_{\RR} u_{\RR} k_1$. Combining the two terms together, and using the equation satisfied by the derivatives on $\Od$, we find

\[
    \frac{\partial \mathcal{K}(\Od_1,k_1)}{\partial k_1} = 2\int_{\Od_1} (\rho_{\DD} u_{\DD} - \rho_{\RR} u_{\RR}) k_1.
\]

\end{proof}

\subsection{Proofs for the case with uncertainties}

\begin{proof}[Proof of the expression of the expected value of the functional: proposition \ref{prop:expectedvalue}]

The computation is done through the Fubini theorem (all the quantities are $\mathrm{L}^2$ and positive) after developing the term and arranging them to factor out the powers of $Y$. Indeed, let us keep in mind that $\bE(Y_i) =0$, $\bE(Y_iY_j) = \delta_{i,j}$, $\bE(Y_iY_jY_kY_l) = 1$ when the indices pair up and is equal to $\bE(Y^4)$ if they are all equal. Otherwise, it is $0$.

\begin{multline*}
    \bE(\mathcal{E}_{\mathcal{K}}(\Od_1,k_1)) 
    = \frac{1}{2}\int_{\Omega}\int_{\mathcal{D}} \left| \sum_{m = 0}^M (\nabla u_{\RR_m} - \nabla u_{\DD_m}) Y_m(\omega) \right|^2 \, \\
    = \int_{\Omega}\int_{\mathcal{D}} \bigg( 
    \sum_{m=0}^M |\nabla u_{\DD_m} - \nabla u_{\RR_m}|^2 Y_m(\omega)^2  \\ 
    + \sum_{0 \le m \ne l \le M} (\nabla u_{\DD_l} - \nabla u_{\RR_l}) \cdot \overline{(\nabla u_{\DD_m} - \nabla u_{\RR_m})} Y_l(\omega) Y_m(\omega)
    \bigg)\, \\
    = \frac{1}{2} \int_{\mathcal{D}} \sum_{m=0}^M |\nabla u_{\DD_m} - \nabla u_{\RR_m}|^2 \,.
\end{multline*}

We recognize the sum of $M + 1$ terms of the form $\int_{\D} |\nabla u_{\DD} - \nabla u_{\RR}|^2$. Using the results of the deterministic case (propositions~\ref{prop:derivativefirststrategy} and \ref{prop:derivativekfirststrategy}), we deduce the desired formula.

\end{proof}
\begin{proof}[Proof of the expression of the variance of the functional: proposition~\ref{prop:variance}]

We proceed in the same way as in the previous proof. The expression is more complex, but still entirely computable numerically.

\begin{align*}
    \bE\left( \mathcal{E}_{\mathcal{K}}(\Od,\omega)^2 \right) 
    &= \frac{1}{4} \int_{\Omega} \left( \int_{\mathcal{D}} 
    \left| \sum_{m=0}^M (\nabla u_{\RR_m} - \nabla u_{\DD_m}) Y_m(\omega) \right|^2 \right)^2.
\end{align*}

  We use Fubini's theorem and the particular case where~\( Y \sim \mathcal{N}(0,1) \): we use $\mathbb{E}(Y_m) = 0$, \( \bE(Y_m^3) = 0 \), \( \bE(Y_m^4) = 3 \).

\begin{multline*}
    \bE\left( \mathcal{E}_{\mathcal{K}}(\Od,\omega)^2 \right) =
    \frac{1}{4} \int_{\Omega} \left( \int_\D |\nabla \mathrm{d}_0|^2 + 2 \mathfrak{Re} \left(\nabla \mathrm{d}_0 \cdot \sum_{m = 1}^M \overline{\nabla \mathrm{d}_m} Y_m \right) + \left| \sum_{m=1}^M \nabla \mathrm{d}_m Y_m\right|^2\right)^2 \\
    = \frac{1}{4} \left(\int_{\D} |\nabla \mathrm{d}_0|^2\right)^2 + \frac{1}{2} \left(\int_{\D} |\nabla \mathrm{d}_0|^2\right) \int_{\Omega}  \left(\int_{\D} \left| \sum_{m = 1}^M \nabla \mathrm{d}_m Y_m\right|^2\right) \\
    + \int_{\Omega} \left(\int_D \mathfrak{Re} \left(\nabla \mathrm{d}_0 \cdot \sum_{m = 1}^M \overline{\nabla \mathrm{d}_m} Y_m \right) \right)^2 + \frac{1}{4} \int_{\Omega} \left( \int_{\D} \left| \sum_{m=1}^M \nabla \mathrm{d}_m Y_m\right|^2\right)^2
    \\
    = \frac{1}{4} \left(\int_{\D} |\nabla \mathrm{d}_0|^2\right)^2 + \frac{1}{2} \left(\int_{\D} |\nabla \mathrm{d}_0|^2\right)  \left( \sum_{m = 1}^M \int_{\D}|\nabla \mathrm{d}_m|^2\right) \\
    + \sum_{m = 1}^M \left( \int_{\D} \mathfrak{Re}\left( \nabla \mathrm{d}_0 \cdot \nabla \overline{\mathrm{d}_m}\right)\right)^2 + \mathfrak{Re}\left( \sum_{1 \le m < l \le M} \left(\int_{\D} \nabla \mathrm{d}_m \overline{\nabla \mathrm{d}_l}\right)^2 \right) \\
    + \frac{3}{4} \sum_{m = 1}^M \left(\int_{\D} |\nabla \mathrm{d}_m|^2\right)^2 + \frac{1}{2} \sum_{1 \le m < l \le M} \left(\int_{\D} |\nabla \mathrm{d}_m|^2 \right)\left(\int_{\D} |\nabla \mathrm{d}_l|^2 \right)
\end{multline*} With the same notations, we have 

\[ \mathcal{E}_{\mathcal{K}}(\Od_1,k_1)^2 = \frac{1}{4} \left(\sum_{m = 0}^M \int_{\D} |\nabla \mathrm{d}_m|^2 \right)^2 = \frac{1}{4} \sum_{m = 0}^M \left( \int_{\D} |\nabla \mathrm{d}_m|^2 \right)^2 + \frac{1}{2} \sum_{1 \le m < l \le M} \left(\int_{\D} |\nabla \mathrm{d}_m|^2 \right)\left(\int_{\D} |\nabla \mathrm{d}_l|^2\right). \] 

Thus, 
\begin{multline*}
    \mathbb{V} \left( \mathcal{E}_{\mathcal{K}}(\Od,\omega) \right) = \frac{1}{4} \left(\int_{\D} |\nabla \mathrm{d}_0|^2\right)^2 + \frac{1}{2} \left(\int_{\D} |\nabla \mathrm{d}_0|^2\right)  \left( \sum_{m = 1}^M \int_{\D}|\nabla \mathrm{d}_m|^2\right) \\
    + \sum_{m = 1}^M \left( \int_{\D} \mathfrak{Re}\left( \nabla \mathrm{d}_0 \cdot \nabla \overline{\mathrm{d}_m}\right)\right)^2 + \mathfrak{Re}\left( \sum_{1 \le m < l \le M} \left(\int_{\D} \nabla \mathrm{d}_m \overline{\nabla \mathrm{d}_l}\right)^2 \right) \\
    + \frac{1}{2} \sum_{m = 1}^M \left(\int_{\D} |\nabla \mathrm{d}_m|^2\right)^2
\end{multline*}

First, we denote $u_{\DD_i}'$ and $u_{\RR_i}'$ the derivatives of the function $\Od_1 \mapsto u_{\DD_i}(\Od_1)$ and $\Od_1 \mapsto u_{\RR_i}(\Od_1)$ as computed in the deterministic case. The expression becomes 

\begin{multline*}
    \mathrm{D}_{\Od_1} \mathbb{V} \left((\Od_1,k_1) \right)(\theta) = \frac{1}{2} \left(\int_{\D} |\nabla \mathrm{d}_0|^2\right)\mathrm{D}_{\Od_1} \mathcal{E}_\mathcal{K}(\Od_1,k_1)(\theta) + 2\sum_{m = 1}^M \left( \int_{\D} \mathfrak{Re}\left( \nabla \mathrm{d}_0 \cdot \nabla \overline{\mathrm{d}_m}\right)\right)\left(\int_{\D} \mathfrak{Re}\left(\nabla \mathrm{d}_0 \cdot \overline{\nabla \mathrm{d}_m'} \right)\right) \\
    + 2 \sum_{1 \le m < l \le M} \mathfrak{Re}\left(\int_{\D} \nabla \mathrm{d}_m \overline{\nabla \mathrm{d}_l}\right) \mathfrak{Re}\left(\int_{D} \nabla \mathrm{d}_m \overline{\nabla \mathrm{d}_l'} + \nabla \mathrm{d}_m' \overline{\nabla \mathrm{d}_l}\right) \\
    + \sum_{m = 1}^M \left(\int_{\D} |\nabla \mathrm{d}_m|^2 \right)\mathfrak{Re}\left( \int_{\partial \Od_1} (k_1^2 - k_2^2) \mathrm{p}_{m,m} (\theta \cdot n) \right). 
\end{multline*}
    
Finally, the expression takes the form \[\mathrm{D}_{\Od_1} \mathbb{V} \left((\Od_1,k_1) \right)(\theta) = \mathfrak{Re}\left( \int_{\partial \Od_1} (k_1^2-k_2^2) \mathcal{G} (\theta \cdot n)\right),\] where 
\begin{multline*}
    \mathcal{G} = \frac{1}{2} \left(\int_{\D} |\nabla \mathrm{d}_0|^2\right)\sum_{m = 1}^M \mathrm{p}_{m,m} + 2\sum_{m = 1}^M \left( \int_{\D} \mathfrak{Re}\left( \nabla \mathrm{d}_0 \cdot \nabla \overline{\mathrm{d}_m}\right)\right)\mathrm{p}_{0,m} \\
    + 4 \sum_{1 \le m < l \le M} \mathfrak{Re}\left(\int_{\D} \nabla \mathrm{d}_m \overline{\nabla \mathrm{d}_l}\right)  \mathrm{p}_{m,l} + \sum_{m = 1}^M \left(\int_{\D} |\nabla \mathrm{d}_m|^2 \right) \mathrm{p}_{m,m}.
\end{multline*}

We proceed in the same way to derivate with respect to the wavenumber $k$.
\end{proof}

\medskip
\paragraph{Acknowledgement.}
F. Caubet was funded by the ANR project STOIQUES (ANR-24-CE40-2216) and the PHC TOUBKAL project 51628SC.


\medskip

\begin{thebibliography}{10}

\bibitem{Afraites2007}
L.~Afraites, M.~Dambrine, K.~Eppler, and D.~Kateb.
\newblock Detecting perfectly insulated obstacles by shape optimization techniques of order two.
\newblock {\em Discrete and Continuous Dynamical Systems - B}, 8(2):389--416, 2007.

\bibitem{allaire2007}
G.~Allaire.
\newblock {\em Conception optimale de structures}, volume~58 of {\em Math{\'e}matiques \& Applications}.
\newblock Springer, Paris, 2007.

\bibitem{Audibert2022}
L.~Audibert, H.~Haddar, and X.~Liu.
\newblock An accelerated level-set method for inverse scattering problems.
\newblock {\em SIAM Journal on Imaging Sciences}, 15(3):1576--1600, 2022.

\bibitem{bonnafont2023inverse}
T.~Bonnafont and F.~Caubet.
\newblock Inverse scattering using a kohn-vogelius formulation and shape optimization method.
\newblock In {\em 2023 17th European Conference on Antennas and Propagation (EuCAP)}, pages 1--5. IEEE, 2023.

\bibitem{Bourgeois2012}
L.~Bourgeois, N.~Chaulet, and H.~Haddar.
\newblock On simultaneous identification of the shape and generalized impedance boundary condition in obstacle scattering.
\newblock {\em SIAM Journal on Scientific Computing}, 34(3):A1824--A1848, 2012.

\bibitem{Boutarene2024}
K.~E.-G. Boutarene, S.~Galleze, and V.~P{\'e}ron.
\newblock {On Ventcel-type transmission conditions for a Helmholtz problem with a non-uniform thin layer}.
\newblock {\em {SeMA Journal: Boletin de la Sociedad Espa{\~n}ola de Matem{\'a}tica Aplicada}}, 2024.

\bibitem{bucci1997electromagnetic}
O.~Bucci and T.~Isernia.
\newblock Electromagnetic inverse scattering: Retrievable information and measurement strategies.
\newblock {\em Radio science}, 32(6):2123--2137, 1997.

\bibitem{cakoni2022inverse}
F.~Cakoni, D.~Colton, and H.~Haddar.
\newblock {\em Inverse scattering theory and transmission eigenvalues}.
\newblock SIAM, 2022.

\bibitem{cakoni2011linear}
F.~Cakoni, D.~Colton, and P.~Monk.
\newblock {\em The linear sampling method in inverse electromagnetic scattering}.
\newblock SIAM, 2011.

\bibitem{cakoni2014qualitative}
F.~Cakoni and D.~L. Colton.
\newblock {\em A qualitative approach to inverse scattering theory}, volume 767.
\newblock Springer, 2014.

\bibitem{Caubet2013}
F.~Caubet, M.~Dambrine, D.~Kateb, and C.~Z. Timimoun.
\newblock {A Kohn-Vogelius formulation to detect an obstacle immersed in a fluid}.
\newblock {\em {Inverse Problems and Imaging }}, 7(1):123--157, 2013.

\bibitem{Chaabane2023}
S.~Chaabane, H.~Haddar, and R.~Jerbi.
\newblock {A combination of Kohn-Vogelius and DDM methods for a geometrical inverse problem}.
\newblock {\em {Inverse Problems}}, 39(9):095001, July 2023.

\bibitem{colton2013}
D.~Colton and R.~Kress.
\newblock {\em Integral Equation Methods in Scattering Theory}.
\newblock Society for Industrial and Applied Mathematics, Philadelphia, PA, 2013.

\bibitem{coltonkress2019}
D.~Colton and R.~Kress.
\newblock {\em Inverse Acoustic and Electromagnetic Scattering Theory}, volume~93 of {\em Applied Mathematical Sciences}.
\newblock Springer, New York, 4 edition, 2019.

\bibitem{crocco2012linear}
L.~Crocco, I.~Catapano, L.~Di~Donato, and T.~Isernia.
\newblock The linear sampling method as a way to quantitative inverse scattering.
\newblock {\em IEEE Transactions on Antennas and Propagation}, 60(4):1844--1853, 2012.

\bibitem{Dambrine2022}
M.~Dambrine and V.~Karnaev.
\newblock {Robust obstacle reconstruction in an elastic medium}.
\newblock {\em {Discrete and Continuous Dynamical Systems - Series B}}, 29(1):151--173, 2024.

\bibitem{Darde2010}
J.~Dard{\'e}.
\newblock {\em {M{\'e}thodes de quasi-r{\'e}versibilit{\'e} et de lignes de niveau appliqu{\'e}es aux probl{\`e}mes inverses elliptiques.}}
\newblock Theses, {Universit{\'e} Paris-Diderot - Paris VII}, Dec. 2010.

\bibitem{Darde2013}
J.~Dard\'e, N.~Hyv\"onen, A.~Sepp\"anen, and S.~Staboulis.
\newblock Simultaneous recovery of admittivity and body shape in electrical impedance tomography: an experimental evaluation.
\newblock {\em Inverse Problems}, 29(8):085004, 16, 2013.

\bibitem{Dorn2000}
O.~Dorn, E.~L. Miller, and C.~M. Rappaport.
\newblock A shape reconstruction method for electromagnetic tomography using adjoint fields and level sets.
\newblock {\em Inverse Problems}, 16(5):1119, oct 2000.

\bibitem{ErnGuermond2021}
A.~Ern and J.-L. Guermond.
\newblock {\em Finite Elements II: Galerkin Approximation, Elliptic and Mixed PDEs}.
\newblock Springer, 2021.

\bibitem{Faverjon2017}
B.~Faverjon, B.~Puig, and T.~Baranger.
\newblock Identification of boundary conditions by solving cauchy problem in linear elasticity with material uncertainties.
\newblock {\em Computers \& Mathematics with Applications}, 73(3):494--504, 2017.

\bibitem{Feppon2020}
F.~Feppon.
\newblock {\em Shape and topology optimization of multiphysics systems}.
\newblock PhD thesis, Thèse de doctorat de l'Université Paris-Saclay préparée à l'école polytechnique, 2019.

\bibitem{Furuya2022}
T.~Furuya and R.~Potthast.
\newblock Inverse medium scattering problems with kalman filter techniques.
\newblock {\em Inverse Problems}, 38(9):095003, aug 2022.

\bibitem{hahner}
P.~H\"ahner.
\newblock A periodic {F}addeev-type solution operator.
\newblock {\em J. Differential Equations}, 128(1):300--308, 1996.

\bibitem{Hecht2012}
F.~Hecht.
\newblock New development in freefem++.
\newblock {\em J. Numer. Math.}, 20(3-4):251--265, 2012.

\bibitem{henrotpierre2018}
A.~Henrot and M.~Pierre.
\newblock {\em Shape Variation and Optimization}, volume~28 of {\em EMS Tracts in Mathematics}.
\newblock European Mathematical Society, Z{\"u}rich, 2018.

\bibitem{kress2018}
R.~Kress.
\newblock Inverse scattering from an impedance boundary condition.
\newblock {\em Mathematical Methods in the Applied Sciences}, 41(16):6343--6354, 2018.

\bibitem{Litman1998}
A.~Litman, D.~Lesselier, and F.~Santosa.
\newblock Reconstruction of a two-dimensional binary obstacle by controlled evolution of a level-set.
\newblock {\em Inverse Problems}, 14(3):685, jun 1998.

\bibitem{Nachman}
A.~I. Nachman.
\newblock Reconstructions from boundary measurements.
\newblock {\em Ann. of Math. (2)}, 128(3):531--576, 1988.

\bibitem{Nesterov1983}
Nesterov and Y.~E.
\newblock A method for solving the convex programming problem with convergence rate {$O(1/k\sp{2})$}.
\newblock {\em Dokl. Akad. Nauk SSSR}, 269(3):543--547, 1983.

\bibitem{nesterov}
Y.~E. Nesterov.
\newblock A method for solving the convex programming problem with convergence rate {$O(1/k\sp{2})$}.
\newblock {\em Dokl. Akad. Nauk SSSR}, 269(3):543--547, 1983.

\bibitem{Novikov}
R.~G. Novikov.
\newblock A multidimensional inverse spectral problem for the equation {$-\Delta\psi +(v(x)-Eu(x))\psi=0$}.
\newblock {\em Funktsional. Anal. i Prilozhen.}, 22(4):11--22, 96, 1988.

\bibitem{Donoghue2015}
O'Donoghue, Brendan, and E.~Cand\`es.
\newblock Adaptive restart for accelerated gradient schemes.
\newblock {\em Found. Comput. Math.}, 15(3):715--732, 2015.

\bibitem{sylvesteruhlmann}
J.~Sylvester and G.~Uhlmann.
\newblock A global uniqueness theorem for an inverse boundary value problem.
\newblock {\em Ann. of Math. (2)}, 125(1):153--169, 1987.

\bibitem{wen2013algorithm}
Y.~Wen, N.~de~Beaucoudrey, J.~Chauveau, and P.~Pouliguen.
\newblock Algorithm for profile function calculation of 3d objects: application for radar target identification in low frequency.
\newblock {\em Progress In Electromagnetics Research B}, 50:273--293, 2013.

\end{thebibliography}
\end{document}